\documentclass[a4paper,11pt]{amsart}
\usepackage[applemac]{inputenc}

 \DeclareMathOperator{\supp}{supp}
\newcommand{\R}{\mathbb{R}}

\newcommand{\N}{\mathbb{N}}
\newtheorem{theorem}{Theorem}
\newtheorem{proposition}[theorem]{Proposition}
\newtheorem{lemma}[theorem]{Lemma}

\theoremstyle{remark}
\newtheorem{remark}{Remark}
\theoremstyle{remark}
\newtheorem{notation}{Notation}

\def\div{\, \mbox{div}\,  }    \def \e {\varepsilon} 
\def\tilde{\widetilde} 
\def\hat{\widehat}
\def\dq{\Delta_q} 
 
\def\d{\partial} 
\def\ov{\overline} 
 
\def\cP{{\mathcal P}}

\def\cC{{\mathcal C}}

\def\div{\, \mbox{\rm div}\,  }

\def\Supp{\, \mbox{\rm Supp}\,  } 
\newenvironment{p}{ 
\noindent\textit{\textbf{Proof:}}~} 
{\hfill\rule{2.1mm}{2.1mm} 
\smallbreak} 
\newcommand{\Sum}{\displaystyle \sum} 
 
\newcommand{\Int}{\displaystyle \int} 
\newcommand{\Frac}{\displaystyle \frac}

\newcommand{\du}{\delta\! u} 
 
\newcommand{\dU}{\delta\! U} 
\newcommand{\dt}{\delta\!\theta} 
 
\newcommand\dPi{\delta\!\Pi} 
 
\def \epsilon {\varepsilon}        
\begin{document}
\title{Global existence results for the anisotropic
 Boussinesq system in dimension two}
\author{Rapha\"el Danchin$^1$ and Marius Paicu$^2$}
\thanks{$^1$Universit\'e Paris-Est, Laboratoire d'Analyse 
et de Math\'ematiques Appliqu\'ees, UMR 8050,
 61 avenue du G\'en\'eral de Gaulle,
94010 Cr\'eteil Cedex, France.  E-mail: danchin$@$univ-paris12.fr}
\thanks{$^2$Universit\'e Paris 11, Laboratoire de Math\'ematiques,  
 B\^atiment 425, 91405 Orsay Cedex, France. 
 E-mail: marius.paicu$@$math.u-psud.fr}

\date\today

\begin{abstract}
We study the global existence issue 
for the two-dimensional Boussinesq system with
\emph{horizontal} viscosity in only one equation.
We first examine the case where 
the Navier-Stokes equation \emph{with no vertical viscosity}
is coupled with a transport equation. 
Second, we consider a coupling between the
classical two-dimensional incompressible Euler equation
and a transport-diffusion equation with 
diffusion \emph{in the horizontal direction only.} 
For the both systems and for arbitrarily large data,  we construct global weak solutions
\emph{\`a la Leray}. Next, we state global well-posedness results for more
regular data. Our results strongly rely on the fact
that the diffusion occurs in a direction perpendicular to the buoyancy force.
 \end{abstract} 
\maketitle 

\section{Introduction}

The Boussinesq system describes the influence of the convection (or convection-diffusion) phenomenon in a viscous or inviscid  fluid. 
It is used as  a toy model for geophysical fluids whenever 
rotation and stratification play an important role  (see for example  J. Pedlosky's book
\cite{Ped}). 
  In the two-dimensional case, the Boussinesq system reads: 
$$
\begin{cases}
\partial_t\theta+u\cdot\nabla \theta-\kappa \Delta\theta=0\\
\partial_t u+u\cdot\nabla u-\nu\Delta u+\nabla\Pi=\theta\, e_2
\quad\hbox{with}\quad e_2=(0,1),\\
\div u=0.
\end{cases}\leqno(B_{\kappa,\nu})
$$
Above, $u=u(t,x)$ denotes the velocity vector-field and  $\theta=\theta(t,x)$
 is a scalar quantity  such 
as  the concentration   of a chemical substance or  the temperature variation in a 
gravity field, in which case $\theta\,e_2$ represents the buoyancy force.
 The nonnegative parameters  $\kappa$ and  $\nu$ denote respectively the molecular diffusion and the viscosity.
In order to simplify the presentation, we restrict ourselves to the whole
plan  case  (that is the space variable $x$ describes the whole $\R^2$) and
focus on the  evolution for positive times (i.e. $t\in\R_+$).  

In the case where both  $\kappa$ and $\nu$ are positive,  classical methods allow to establish the global existence of regular solutions (see for example  \cite{CD,Guo}).
On the other hand, if
$\kappa=\nu=0$ then  constructing  global unique solutions  for some  
nonconstant $\theta_0$  
is  a challenging open problem (even in the two-dimensional case) 
which has many similarities with the global existence problem for
the three-dimensional incompressible Euler equations.

The intermediate situation where the diffusion acts only on one of 
the equations  has been investigated in a number of recent papers. 
Under various regularity assumptions on the initial data,  it has been shown that
for arbitrarily large initial data,
systems $(B_{\kappa,0})$ with $\kappa>0$ and $(B_{0,\nu})$ with $\nu>0$ 
admit a global unique solution
 (see for example  \cite{AH,dingo,DP1,DP2,DP3,HK1,HK2}).
\smallbreak
In the present paper, we aim at making one more step toward the study 
of the system with $\kappa=\nu=0$
by assuming that the diffusion or the viscosity occurs in the
horizontal direction and in one of the equations only. 
 More precisely, we want to consider the following two systems:
\begin{equation}\label{eq:b-vitesse}
\begin{cases}
\partial_t\theta+u\cdot\nabla \theta=0\\
\partial_t u+u\cdot\nabla u-\nu \partial_{1}^2 u+\nabla\Pi=\theta\, e_2\\
\div u=0
\end{cases}
\end{equation}
and
\begin{equation}\label{eq:b-temperature}
\begin{cases}
\partial_t\theta+u\cdot\nabla\theta-\kappa\partial_{1}^2\theta=0\\
\partial_t u+u\cdot\nabla u+\nabla\Pi=\theta e_2\\
\div u=0.
\end{cases}
\end{equation}
Let us stress that the anisotropic viscosity or diffusion assumptions
are consistent with the study of geophysical fluids. 
It turns out that, in certain regimes and after suitable rescaling, 
the vertical viscosity (or diffusion) is negligible with respect to the horizontal viscosity
(or diffusion) (for more details, one may refer to  \cite{CDGG06}). 
For the standard three-dimensional incompressible  Navier-Stokes equations, 
the first mathematical results concerning anisotropic viscosity  have been obtained in 
\cite{Chemin,Paicu}. 
\smallbreak
On the one hand, it is clear that small variations over the classical methods for solving quasi-linear hyperbolic
systems would give local well-posedness for Systems \eqref{eq:b-vitesse}
and  \eqref{eq:b-temperature} with initial data in Sobolev spaces 
with suitably large index. 
On the other hand, since diffusion occurs in only one direction and one
equation, it is not obvious that those solutions are actually global. 
The present paper is dedicated to the study of global existence for the initial 
value problem associated to Systems 
 \eqref{eq:b-vitesse}
and  \eqref{eq:b-temperature} with (possibly) large initial data. 

In order to state our main result pertaining to System  \eqref{eq:b-vitesse}, let
us introduce the set $\sqrt L$ of those functions 
$f$ which  belong to every  space  $L^p$ with  $2\leq p<\infty$  and satisfy
\begin{equation}\label{eq:racineL}
\|f\|_{\sqrt L}:=\sup_{p\geq2}\,\frac{\|f\|_{L^p}}{\sqrt{p\!-\!1}}<\infty.
\end{equation}
\begin{theorem}\label{th:resultat1}
Let $s\in]1/2,1].$  For all function   $\theta_0\in H^s\cap L^\infty$ and divergence free
vector-field $u_0\in H^1$ with vorticity  $\omega_0:=\d_1v_0^2-\d_2v_0^1$ in  $\sqrt L,$
System~$\eqref{eq:b-vitesse}$ with data $(\theta_0,u_0)$ admits a unique global solution 
 $(\theta,u)$ such that\footnote{In all the paper, 
we agree that if  $X$ is a reflexive  Banach space,
and  $I\subset\R,$ an  interval then
 $\cC_w(I;X)$ stands
for the set of  weakly continuous  functions on  $I$ with values in $X.$}
$\theta\in\cC_w(\R_+; L^\infty)\cap\cC(\R_+;H^{s-\e})$ for all  $\e>0$  and
$$u\in\cC_w(\R_+;H^1),\quad \omega\in L^\infty_{loc}(\R_+;\sqrt L)\ \hbox{ and }\ 
\nabla u\in L^2_{loc}(\R_+;\sqrt L).$$
\end{theorem}
\begin{remark}
The assumption that  $\theta_0\in H^s(\R^2)$ for some   $s>\frac 12$ 
is needed for uniqueness only. 
It turns out that for less  regular  
 initial data one can construct  
finite energy global weak solutions to System  \eqref{eq:b-vitesse}, 
in the spirit of those which have been obtained by J. Leray for 
the standard Navier-Stokes equation in his pioneering paper \cite{Leray}.

We shall also establish a global well-posedness results
for smooth initial data.
 \end{remark}
 Let us now state our main result pertaining 
 to System \eqref{eq:b-temperature}:
\begin{theorem}\label{th:resultat2} 
Let $1<s<\frac32$ and  $\theta_0\in H^1$ such that 
 $|\partial_1|^{s} \theta_0\in L^2.$  Let  $u_0$ be a divergence free vector-field
 with coefficients in  $H^1$ and vorticity  $\omega_0$ in $L^\infty.$
Then System $\eqref{eq:b-temperature}$ with initial data $(\theta_0,u_0)$ 
admits a global  unique solution $(\theta,u)$
in $\cC_w(\R_+;H^1)$ such that, in addition, 
 $$
 \displaylines{
\theta\in L^\infty(\R_+;H^1), \quad 
|\partial_1|^s\theta\in L^\infty(\R_+; L^2), \quad\omega\in L^\infty_{loc}(\R_+;L^\infty)\cr
\partial_1\theta\in L^2(\R_+;H^1),\quad 
|\partial_1|^{1+s}\theta\in L^2_{loc}(\R_+;L^2).}
$$
\end{theorem}
\begin{notation}
In the above statement, operator $|\d_1|^\sigma$ is
defined as follows:
$$
|\d_1|^\sigma f(x):=\frac1{4\pi^2}\int e^{ix\cdot\xi} |\xi_1|^\sigma \widehat f(\xi)\,d\xi.
$$
\end{notation}
\begin{remark}
Further in the paper, we shall also state the global existence 
of finite energy weak solutions corresponding to  less regular initial data  (see Theorem \ref{th:resultat4}).
\end{remark}
\begin{remark}
Like in \cite{DP3}, in all the statements pertaining to System \eqref{eq:b-temperature},
 one may replace the assumption that  $u_0\in L^2$
(which is slightly restrictive since in the case $\omega_0\in L^1$ 
it implies that $\int_{\R^2}\omega_0\,dx=0$)
by $(u_0-\sigma)\in L^2$ for some fixed smooth stationary 
radial solution  $\sigma$ to Euler equations. 
For the sake of simplicity however, we here restrict ourselves to the case
where  $u_0$ is in~$L^2.$
\end{remark}
Let us emphasize that, in contrast with the previous studies devoted
to the Boussinesq system,  all the results presented here 
strongly rely on the fact that the buoyancy force is \emph{vertical}. 
As a matter of fact, it is not clear at all that if the direction 
of the force is changed then having horizontal 
diffusion allows to  improve significantly the results compared 
to the case $\kappa=\nu=0.$

Let us now briefly explain where the direction of the buoyancy
force comes into play, and give an insight of the main
arguments that we used in the proofs. 
In the case of System \eqref{eq:b-vitesse}, the vorticity 
$\omega:=\d_1u^2-\d_2u^1$ satisfies
$$
\d_t\omega+u\cdot\nabla\omega-\nu\d_1^2\omega=\d_1\theta,
$$
from which we  easily get (at least formally) 
$$
\frac12\frac d{dt}\|\omega\|_{L^2}^2+\nu\|\d_1\omega\|_{L^2}^2
=-\int_{\R^2}\theta \: \partial_1\omega\; dx.
$$
Taking advantage of the Young inequality and of the fact that
$\|\theta(t)\|_{L^2}\leq\|\theta_0\|_{L^2},$ it
is thus possible to get a bound for $\omega$ in $L^\infty_{loc}(\R_+;L^2).$
In fact, it turns out that similar arguments enable us 
to bound stronger norms of the solution so that
it will be possible to prove the global existence part of Theorem \ref{th:resultat1}.
As the velocity field that we have constructed 
fails to be Lipschitz, proving uniqueness requires our
using  losing estimates
for transport or transport-diffusion equations in the spirit of  \cite{BCh,DANCHI,DP1,DP2}. 
\smallbreak
If we consider System \eqref{eq:b-temperature} then
the vorticity equation reduces to 
$$
\d_t\omega+u\cdot\nabla\omega=\d_1\theta,
$$
so that one may write
$$
\frac12\frac d{dt}\|\omega\|_{L^2}^2
=\int \partial_1\theta\: \omega \;dx.
$$
Now, it turns out that the temperature equation in System
\eqref{eq:b-temperature} provides us with the following  bound:  
$$\|\theta(t)\|_{L^2}^2+
2\kappa\int_0^t\|\partial_1\theta(\tau)\|^2_{L^2(\R^2)}d\tau\leq \|\theta_0\|_{L^2}^2.$$
Therefore, using   Young's inequality, it is still
possible to get a global control on $\omega$ in $L^\infty_{loc}(\R_+;L^2).$
In order to get a global result \emph{with uniqueness} however, 
we have to consider initial data with much more regularity. 
Indeed, we have to observe that System
\eqref{eq:b-temperature} contains the Euler system as a particular case
(just take $\theta\equiv0$) so that, according to Yudovich result in \cite{Yudo}  one can 
barely expect 
to have  uniqueness if the vorticity is not in $L^\infty.$ 
Now, from the vorticity equation, it is clear that bounding the vorticity
in $L^\infty$ requires that 
$\partial_1\theta$ is in $L^1_{loc}(\R_+;L^\infty)$. If we assume that $\theta_0\in H^1)$
then  we shall be able to prove that 
the horizontal smoothing effect ensures that $\nabla\theta$ is in $L^2_{loc}(\R_+;H^1).$
As the space   $H^1$ fails to be embedded in $L^\infty$ however, 
 we will have to assume even more regularity in the horizontal direction, 
and to check that this additional regularity is preserved for all positive time.
These plain considerations explain  the assumptions made in the statement of Theorem 
\ref{th:resultat2}.
As for uniqueness, it will follow from adaptations of the Yudovich method in  \cite{Yudo}. 
\medbreak
The paper unfolds as follows: section \ref{visco:horizontale:vitesse}
 is devoted to the study of System \eqref{eq:b-vitesse} 
whereas Section \ref{s:b-temperature} deals with System \eqref{eq:b-temperature}. 
A few  technical lemmas have been postponed in 
the final section of the paper (in particular losing a priori estimates
for transport equation with anisotropic diffusion). 

As usual, we agree that  $C$ denotes a harmless positive
constant, the meaning of which is clear from the context.


\section{The case of an horizontal viscosity} \label{visco:horizontale:vitesse}

This part is devoted to the study of the initial value problem
for System \eqref{eq:b-vitesse} under   various
regularity hypotheses. We aim at getting global
results for possibly large data.

More precisely, in the first subsection, we prove that for
any data  $(\theta_0,u_0)\in L^2\times H^1$ with  $\div u_0=0,$ 
System \eqref{eq:b-vitesse} admits at least one global solution 
with finite energy. 
The next subsection is devoted to a local well-posedness result
for smooth data, together
with a continuation criterion involving the $L^\infty$ norm
of  $(\theta,\nabla u).$ 
In subsection \ref{ss:srg}, we state a sharper continuation criterion
involving a \emph{weaker} norm
of the velocity which  is (formally) controlled for
all time by System \eqref{eq:b-vitesse}. 
This will enable us to state that the system is globally 
well-posed in Sobolev spaces with  large enough index
(see  Theorem \ref{th:resultat3}). 
The last two subsections are devoted to the proof
of our global existence and uniqueness result for rough data
(namely Theorem \ref{th:resultat1}).

\subsection{Global weak solutions}

In order to motivate our statement, let
us first write out the ``natural'' energy estimates
associated to System
\eqref{eq:b-vitesse}.

On the one hand, because  $\div u=0,$ we have
\begin{equation}\label{eq:theta}
\|\theta(t)\|_{L^2}\leq\|\theta_0\|_{L^2}.
\end{equation} 
On the other hand, 
taking the  $L^2$ inner product of the velocity
equation with $u,$ we find that
$$
\frac12\frac
d{dt}\biggl(\|u(t)\|_{L^2}^2+2\nu\int_0^t\|\partial_1u\|_{L^2}^2\,d\tau\biggr)
\leq\|\theta\|_{L^2}\|u\|_{L^2}. $$
Using Gronwall lemma and 
\eqref{eq:theta}, we thus get
\begin{equation}\label{eq:energie}
\|u(t)\|_{L^2}^2+2\nu\int_0^t\|\partial_{1}u\|_{L^2}^2\,d\tau\leq 
\bigl(\|u_0\|_{L^2}+t\,\|\theta_0\|_{L^2}\bigr)^2.
\end{equation}
Let us stress the fact that the 
above energy bounds imply that all the components of 
$\nabla u$ except $\d_2u^1$ are smoothed out 
for positive time. 
Indeed, combining the  $L^2_{loc}(\R_+;L^2)$
bound for $\partial_1u$ which is available from 
\eqref{eq:energie}
with the fact that $\div u=0$ ensures that 
 $\d_1u^1$ $\d_1u^2$ and $\d_2u^2$
 are in  $L^2_{loc}(\R_+;L^2).$
 However, as the last component $\d_2u^1$ is unlikely to
 be bounded in  $L^2_{loc}(\R_+;L^2)$ if no stronger assumption, 
it is not clear that 
one may construct global weak solutions for $L^2$ data
(in contrast with the standard Navier-Stokes equations \cite{Leray}
or with the Boussinesq system with isotropic
viscosity  \cite{DP1}).

This induces us to consider initial velocity fields in $H^1.$ Now,
in order to get a global bound for the $H^1$ norm
of the velocity, one may consider  the vorticity
equation
\begin{equation}\label{eq:tourbillon}
\partial_t \omega+u\cdot\nabla\omega-\nu\partial_{1}^2 
\omega=\partial_{1}\theta.
\end{equation}
Combining an energy method
with Young's inequality, we get
 $$\frac
 12\frac{d}{dt}\|\omega\|_{L^2}^2+\nu\|\partial_{1}\omega\|_{L^2}^2
=-\int\theta\,\partial_{1}\omega\,dx
\leq \frac{\nu}{2}\|\partial_{1}\omega\|_{L^2}^2
+\frac{1}{2\nu}\|\theta\|^2_{L^2}.$$
Therefore,
$$\frac{d}{dt}\|\omega\|_{L^2}^2+\nu\|\partial_{1}\omega\|_{L^2}^2
\leq\frac{1}{\nu}\|\theta\|^2_{L^2},$$
whence,  according to \eqref{eq:theta},
\begin{equation}\label{eq:omegaL2}
\|\omega(t)\|_{L^2}^2+\nu\int_0^t
\|\partial_{1}\omega\|_{L^2}^2\,d\tau
\leq\|\omega_0\|_{L^2}^2+\frac{t}{\nu}\|\theta_0\|^2_{L^2}.
\end{equation}
In short, one may formally bound
$u$ in $L^\infty_{loc}(\R_+;H^1)$ 
and  $\theta$ in $L^\infty(\R_+;L^2)$
which leads to the following statement.
\begin{theorem}\label{th:globalfaible}
Let  $u_0\in H^1$ be a divergence free vector-field 
and  $\theta_0\in L^2$. 
System $\eqref{eq:b-vitesse}$ has a unique global solution 
$(\theta,u)$ such that\footnote{We agree that if  $X$ is a  Banach space,
and  $I\subset\R,$ an  interval then
 $\cC_b(I;X)$ stands
for the set of continuous bounded 
   functions on  $I$ with values in $X.$}
\begin{equation}\label{eq:globalfaible}
\theta\in\cC_b(\R_+;L^2),\quad
u\in\cC_w(\R_+;H^1)\quad\hbox{and}\quad
u^2\in L^2_{loc}(\R_+;H^2).
\end{equation}
\end{theorem}
\begin{p}
It is only a matter of making 
the above computations rigorous. 
For that, one may for instance
use the  Friedrichs method:
define the spectral cut-off  $J_n$ by
$$\widehat{J_nf}(\xi)=1_{[0,n]}(|\xi|)\widehat
f(\xi),\qquad n\geq1,$$  
and solve the following ODE
in the space $L^2_n:=\bigl\{f\in(L^2(\R^2))^3\,/\, 
\supp\widehat f\subset B(0,n)\bigr\}$:
\begin{equation}\label{eq:approche}
\begin{cases}
\partial_t\theta+J_n\div(J_n\theta\,J_nu)=0,\\
\partial_tu+\cP J_{n}\div(\cP J_{n}u\otimes \cP J_{n}u)
-\nu\d_1^2\cP J_{n}u=\cP J_{n}(\theta e_2),\\
(\theta,u)|_{t=0}=J_{n}(\theta_0,u_0).
\end{cases}
\end{equation}
{}From the Cauchy-Lipschitz theorem, 
we get a unique  maximal solution $(\theta_n,u_n)$
in $\cC^1([0,T_n^*[;L^2_n).$
Because $J_n^2=J_n,$ $\cP^2=\cP$ and
$ J_n\cP=\cP J_n,$ 
we discover that  $(\theta_n,\cP u_n)$ and  $(J_n\theta_n,J_nu_n)$
are also  solutions.  By uniqueness, 
we thus have  $\cP u_n=u_n$ (i. e. 
 $\div u_n=0$), $J_nu_n=u_n$ and $J_n\theta_n=\theta_n.$
Therefore, 
\begin{equation}\label{eq:approche1}
\begin{cases}
\partial_t\theta_n+J_n\div(\theta_{n}u_n)=0,\\
\partial_tu_{n}+\cP J_{n}\div(u_{n}\otimes u_{n})
-\nu\d_1^2 u_{n}=\cP J_{n}(\theta_{n} e_2),\\
\div u_n=0.
\end{cases}
\end{equation}
As Operators $J_n$ 
and  $\cP J_n$ are orthogonal projectors
for the $L^2$ inner product, 
the above formal calculations remain unchanged.
Therefore, we still have as before
\begin{eqnarray}\label{eq:energie1}
&\|\theta_n(t)\|_{L^2}^2=\|J_n\theta_0\|_{L^2}^2\leq \|\theta_0\|_{L^2}^2,\\[1.5ex]
\label{eq:energie2}
&\|u_n(t)\|_{L^2}^2+2\nu\Int_0^t\|\d_1u_n\|_{L^2}^2\,d\tau
\leq \bigl(\|u_0\|_{L^2}+t\,\|\theta_0\|_{L^2}\bigr)^2.
\end{eqnarray}
This implies that $(\theta_n,u_n)$ remains bounded in $L^2_n$ 
for finite time, whence $T_n^*=+\infty. $

Next, applying the ${\rm curl}\,$ operator
to $\eqref{eq:approche1}_2,$ we get
$$
\d_t\omega_n+J_n(u_n\cdot\nabla\omega_n)
-\nu\d^2_1\omega_n=\d_1\theta_n
\quad\hbox{with}\quad
\omega_n=\d_1u_n^2-\d_2u_n^1.
$$
Arguing as for proving \eqref{eq:omegaL2}, we thus get
\begin{equation}\label{eq:energie3}
\|\omega_n(t)\|_{L^2}^2+\nu\int_0^t\|\partial_1\omega_n\|_{L^2}^2\,d\tau\leq
\|\omega_0\|_{L^2}^2+\frac t\nu\|\theta_0\|_{L^2}^2.
\end{equation}
This implies that $(\omega_n)_{n\in\N}$
is bounded in  $L^\infty_{loc}(\R_+;L^2).$
Now, it is well known 
that the divergence-free property
entails that
$$
\|\nabla u_n\|_{L^2}=\|\omega_n\|_{L^2}\quad\hbox{et}\quad
\Delta u^2_n=\partial_1\omega_n.
$$
So one can conclude that 
\begin{itemize}
\item $(\theta_n)_{n\in\N}$ is bounded in $L^\infty(\R_+;L^2),$
\item $(u_n)_{n\in\N}$ is bounded in 
$L^\infty_{loc}(\R_+;H^1),$
\item $(u_n^2)_{n\in\N}$ is bounded in 
$L^2_{loc}(\R_+;H^2).$
\end{itemize}
This is enough to pass to the limit (up to extraction)
in  \eqref{eq:approche1}.
Indeed, putting together
 the continuous embedding  $H^1\hookrightarrow L^4$
and H\"older inequality,
we see that the first two properties imply
that  $(\theta_nu_n)_{n\in\N}$ is
bounded in $L^4_{loc}(\R_+;L^{\frac43})$ whence, 
by embedding, in $L^4_{loc}(\R_+;H^{-\frac12}).$
Therefore
$(\d_t\theta_n)_{n\in\N}$ 
is bounded in  $L^4_{loc}(\R_+;H^{-\frac32}).$
Likewise, 
$(\d_tu_n)_{n\in\N}$ 
is bounded in  $L^2_{loc}(\R_+;H^{-1}).$
Since the embeddings $H^{-\frac32}\hookrightarrow L^2$ 
and  $H^{-1}\hookrightarrow H^1$ are locally compact, 
the classical  Aubin-Lions argument  (see
 e.g. \cite{Aubin}) allows to conclude that, up to extraction, sequence
$(\theta_n,u_n)_{n\in\N}$ has a limit
$(\theta,u)$ satisfying
System \eqref{eq:b-vitesse} and that
$$
\theta\in L^\infty(\R_+;L^2),\quad u\in L^\infty_{loc}(\R_+;H^1)\ \hbox{ and }\
u^2\in L^2_{loc}(\R_+;H^2).
$$
{}From standard arguments
relying on the time continuity of 
$(\theta,u)$ in low norms, it is easy 
to prove the weak time continuity result.
Finally, since  $\theta$ is transported
by the flow of a divergence free vector-field
with coefficients in $L^2_{loc}(\R_+;H^1),$
we get in addition that $\theta\in\cC(\R_+;L^2)$ (see e.g. 
\cite{DPL}).
\end{p}


\subsection{Local smooth solutions}\label{ss:regloc}
Here we aim at proving the local well-posedness for
System~\eqref{eq:b-vitesse} with  initial data 
$(\theta_0,u_0)$ in  $H^{s-1}\times H^s$ for some  $s>2.$
\smallbreak
The proof will follow from  an energy method once  
the system has been localized in dyadic frequencies. 
This  localization may be done by means of a nonhomogeneous
 Littlewood-Paley decomposition.
In order to define the dyadic blocks $\dq$ used in this decomposition,  one may proceed
as in \cite{C}: starting from a couple $(\chi,\varphi)$
of smooth nonnegative  functions such that
$$\displaylines{\Supp\chi\subset\{\xi\in\R^2\,/\,  |\xi|\leq 4/3\},\qquad
\Supp\varphi\subset\{\xi\in\R^2\,/\,3/4\leq|\xi|\leq8/3\},\cr
\chi(\xi)+\sum_{q\geq0}\varphi(2^{-q}\xi)=1\quad\hbox{for all  }\ 
\xi\in\R^2,}
$$
we set
$$
\Delta_q:=0\ \hbox{ if }\ q\leq-2,\quad
\Delta_{-1}:=\chi(D)\quad\hbox{and}\quad
\Delta_q:=\varphi(2^{-q}D)\ \hbox{ if }\ q\geq0.
$$
We also introduce the low frequency cut-off
$$
S_q:=\sum_{p\leq q-1}\Delta_q
$$
and (for technical purposes) the \emph{modified} low frequency cut-off
$\ov S_q$ defined by 
\begin{equation}\label{eq:Sq}
\ov S_{-1}=\Delta_{-1}=S_0\quad\hbox{and}\quad
\ov S_q=S_q\quad\hbox{if}\quad q\not=-1.
\end{equation}
It may be easily checked that $$
u=\sum_{q\geq-1}\dq u\quad\hbox{for all tempered distribution }\ u
$$
and that the set of tempered distributions $u$ satisfying
$$
\biggl(\sum_{q\geq-1}2^{2qs}\|\dq u\|_{L^2}^2\biggr)^{\frac12}<\infty
$$
coincides with the  Sobolev space $H^s,$  the above 
left-hand side defining a norm equivalent 
to the usual one.
\medbreak
Let us now state the main result of this subsection:
\begin{proposition}\label{p:local}
Let $(\theta_0,u_0)$ be in $H^{s-1}\times H^s$ with  $s>2.$
Assume that  ${\rm div}\, u_0=0.$ There exists a positive time  $T$ depending
only (continuously)  on $\nu$ and on  $\|(\theta_0,\omega_0)\|_{H^{s-1}}$ such that
System $\eqref{eq:b-vitesse}$ admits a unique solution 
$(\theta,u)$ in $\cC([0,T];H^{s-1}\times H^s).$
Moreover, $u^2\in L^2([0,T];H^{s+1}).$
\end{proposition}
\begin{p}
The uniqueness is a straightforward  consequence of 
a more general result  (see Proposition
\ref{p:unicite}) the proof of which is postponed to subsection
\ref{ss:uniqueness}.
So let us focus on the existence part
of the above proposition, which 
is mostly a  consequence  of the
$H^{s-1}\times H^s$ a priori estimates associated to System
\eqref{eq:b-vitesse}.
\subsubsection*{1. A priori estimates in $H^{s-1}\times H^s$ }
Let 
$(\theta,u)\in\cC^1([0,T];H^\infty)$ satisfy  \eqref{eq:b-vitesse}.
We claim that
there exists a constant $C$ depending only on $s$ and such that 
for all $t\in[0,T],$ we have
\begin{equation}\label{eq:apriori}
\|(\theta,\omega)(t)\|_{H^{s-1}}^2
+\nu\int_0^t\|\d_1\omega\|_{H^{s-1}}^2\,d\tau
\leq \|(\theta_0,\omega_0)\|_{H^{s-1}}^2\,
e^{\frac t\nu}\,e^{C\int_0^t\|(\theta,\nabla v)\|_{L^\infty}\,d\tau}.
\end{equation}
Indeed, applying operator
$\Delta_q$ to the equation satisfied by  $\theta$ yields
$$\partial_t\Delta_q\theta+\ov S_{q-1}u\cdot\nabla\Delta_q\theta
=F_q(u,\theta)\quad\hbox{with}\quad
F_q(u,\theta):=\ov S_{q-1}u\cdot\nabla\Delta_q\theta
-\Delta_q(u\cdot\nabla\theta).$$
Taking the  $L^2$ inner product  of the above equality with  $\Delta_q\theta$ 
and using the divergence free condition, we thus get
$$
\frac12\frac d{dt}\|\Delta_q\theta\|_{L^2}^2\leq
\|F_q(u,\theta)\|_{L^2}\|\Delta_q\theta\|_{L^2}.
$$
In the appendix (see Inequality  \eqref{eq:Fqa}), we state that 
$$
\|F_q(u,\theta)\|_{L^2}\leq C\Bigl(\|\nabla u\|_{L^\infty}
\sum_{q'\geq q-4}2^{q-q'}\|\Delta_{q'}\theta\|_{L^2}
+\|\theta\|_{L^\infty}\sum_{|q'-q|\leq1}\|\Delta_{q'}\omega\|_{L^2}
\Bigr).$$
Plugging this inequality in \eqref{eq:theta1} then multiplying both sides by $2^{2q(s-1)}$
and summing up  over $q\geq-1,$ we get
\begin{equation}\label{eq:thetaHs-1}
\frac d{dt}\|\theta\|_{H^{s-1}}^2
\leq C\Bigl(\|\nabla u\|_{L^\infty}\|\theta\|_{H^{s-1}}^2
+\|\theta\|_{L^\infty}\|\omega\|_{H^{s-1}}\|\theta\|_{H^{s-1}}\Bigr).
\end{equation}
In order to get a  $H^{s-1}$ estimate for $\omega,$  
one may apply  $\Delta_q$ to the vorticity equation.
With the above notation, we get
$$
\d_t\dq\omega+\ov S_{q-1}u\cdot\nabla\dq\omega-\nu\d_1\dq\omega=\d_1\dq\theta
+F_q(u,\omega).
$$
Taking the $L^2$ inner product of this inequality with  $\Delta_q\omega$ 
and using once again the divergence free condition, we get
after integration by parts,
$$
\frac12\frac d{dt}\|\Delta_q\omega\|_{L^2}^2+\nu\|\d_1\dq\omega\|_{L^2}^2
=-\int\dq\theta\d_1\dq\omega\,dx
+\int F_q(u,\omega)\dq\omega\,dx.
$$
Now, we notice that, by virtue of the  Young  inequality,
$$
-\!\int\!\dq\theta\d_1\dq\omega\,dx\leq\frac{\|\dq\theta\|_{L^2}^2}{2\nu}
+\frac\nu2\|\d_1\dq\omega\|_{L^2}^2,
$$
and, according to \eqref{eq:Fqc}, 
$$
\|F_q(u,\omega)\|_{L^2}\leq C\|\nabla u\|_{L^\infty}\!
\!\sum_{q'\geq q\!-\!4}2^{q-q'}\|\Delta_{q'}\omega\|_{L^2}.
$$
Therefore
$$
\frac d{dt}\|\Delta_q\omega\|_{L^2}^2+\nu\|\d_1\dq\omega\|_{L^2}^2
\leq \nu^{-1}\|\dq\theta\|_{L^2}^2+C\|\dq\omega\|_{L^2}\|\nabla u\|_{L^\infty}
\!\sum_{q'\geq q-4}2^{q-q'}\|\Delta_{q'}\omega\|_{L^2}.
$$
Multiplying both sides by $2^{2q(s-1)}$ then summing up over $q\geq-1,$ 
we end up with
$$
\frac d{dt}\|\omega\|_{H^{s-1}}^2
+\nu\|\d_1\omega\|_{H^{s-1}}^2
\leq\nu^{-1}\|\theta\|_{H^{s-1}}^2+C\|\nabla u\|_{L^\infty}
\|\omega\|_{H^{s-1}}^2.
$$
It is now clear that adding up this latter
inequality to \eqref{eq:thetaHs-1} then applying Gronwall lemma completes
the proof of Inequality \eqref{eq:apriori}.

\subsubsection*{2. The proof of local existence.}
One can use again the Friedrichs method introduced in
the proof of Theorem \ref{th:globalfaible}. 
As Operators $J_n$ are orthogonal  projectors for all  Sobolev spaces,
they do not modify the energy estimates leading to Inequality \eqref{eq:apriori}.
Therefore, the approximate solution
$(\theta_n,u_n)$ to \eqref{eq:approche1} satisfies 
$$
\|(\theta_n,\omega_n)(t)\|_{H^{s-1}}^2
+\nu\int_0^t\|\d_1\omega_n\|_{H^{s-1}}^2\,d\tau
\leq \|J_n(\theta_0,\omega_0)\|_{H^{s-1}}^2\,
e^{\frac t\nu}\,e^{C\int_0^t\|(\theta_n,\nabla u_n)\|_{L^\infty}\,d\tau}.
$$
Of course, the $L^2$ norm of  $u_n$ is controlled by virtue of  \eqref{eq:energie2}.
As $s>2,$ the space $H^{s-1}$ continuously embeds in  $L^\infty.$
Since  $\|\nabla u_n\|_{H^{s-1}}=\|\omega_n\|_{H^{s-1}},$
 the previous inequality thus entails that
 $$
X_n(t)\leq \|(\theta_0,\omega_0)\|_{H^{s\!-\!1}}\,e^{\frac t{2\nu}}\,
e^{C\int_0^tX_n(\tau)\,d\tau}
\quad\hbox{with}\quad
X_n^2(t):= \|(\theta_n,\omega_n)(t)\|_{H^{s\!-\!1}}^2
+\nu\int_0^t\|\d_1\omega_n\|_{H^{s\!-\!1}}^2\,d\tau.
$$
This inequality may be easily integrated into
$$
\exp\biggl(-C\int_0^tX_n(\tau)\,d\tau\biggr)\geq 1-2C\nu X_0e^{\frac t{2\nu}}
\quad\hbox{for all}\quad t\geq0.
$$
Therefore, there exists a $c>0$ such that if we set 
\begin{equation}\label{eq:temps}
T:=2\nu\log\biggl(\frac c{\nu\|(\theta_0,\omega_0)\|_{H^{s-1}}}\biggr)
\end{equation}
then 
\begin{itemize}
\item $(\theta_n)_{n\in\N}$  is bounded in  $L^\infty([0,T];H^{s-1}),$
\item $(u_n)_{n\in\N}$  is bounded in
$L^\infty([0,T];H^s),$
\item $(u_n^2)_{n\in\N}$  is bounded in 
$L^2([0,T];H^{s+1}).$
\end{itemize}
Mimicking the compactness argument used for proving Theorem
\ref{th:globalfaible}, one can now conclude that there exists a solution 
$(\theta,u)$ satisfying  $(\theta,u)\in L^\infty([0,T];H^{s-1}\times H^s)$
and $u^2\in L^2([0,T];H^{s+1}).$ 
The time continuity follows from the fact that   $\theta$ and $\omega$ 
satisfy  transport equations with  $H^{s-1}$ initial 
data and a $L^2([0,T];H^{s-1})$ source term. 
This completes the proof of Proposition \ref{p:local}.
\end{p}


\subsection{Global smooth solutions}\label{ss:srg}

Here we aim at proving that the local smooth solutions 
which have been constructed so far may 
be extended to all positive time.
Exhibiting a polynomial control of
$\|\nabla u(t)\|_{\sqrt L}$ (where the space $\sqrt L$ has been defined in
 \eqref{eq:racineL})  is the cornerstone of
 this extension. More precisely, we shall 
 first prove that the  $L^1([0,T];LL)$ 
norm of  $\nabla u$ with
\begin{equation}\label{eq:LL}
LL:=\biggl\{f\in {\mathcal S}'\,/\,
\|f\|_{LL}:=\sup_{q\geq0}
\frac{\|S_qf\|_{L^\infty}}{q+1}<\infty\biggr\}
\end{equation}
controls the Sobolev regularity of the  solutions to System \eqref{eq:b-vitesse}.
Next, we shall state that, under the  hypotheses of Proposition \ref{p:local}, 
the  norm of $\nabla u$ in  $L^2([0,T];\sqrt L)$ (which is \emph{stronger} to
the $L^1([0,T];LL)$ norm) may be bounded for all time by a fixed polynomial
the coefficients of which depend only on low norms of the data, and on $\nu.$ 
Combining this with Proposition  \ref{p:local}
will lead to the following global existence statement:
\begin{theorem}\label{th:resultat3}
 Let $(\theta_0,u_0)$ be in  $H^{s-1}\times H^s$ for some $s>2.$
Assume that  ${\rm div}\, u_0=0.$ Then system $\eqref{eq:b-vitesse}$
has a unique global solution $(\theta,u)$ such that
 $$(\theta,u)\in\cC(\R_+;H^{s-1}\times H^s)\quad\hbox{and}\quad
 u^2\in L^2_{loc}(\R_+;H^{s+1}).$$
\end{theorem}
As a first step for proving Theorem \ref{th:resultat3},
let us show the following lemma:
\begin{lemma}\label{l:explosion}
Let $(\theta,u)$ be a  solution to  $\eqref{eq:b-vitesse}$
in $\cC([0,T);H^{s-1}\times H^s)$ with $s>2$.
If 
 $$\int_0^T\|\nabla u(t)\|_{LL}\,dt<\infty$$
then   $(\theta,u)$ may be continued beyond  $T$
into a smooth solution of
 $\eqref{eq:b-vitesse}.$ 
\end{lemma}
\begin{p}
Putting together the lower bound for the lifespan 
of $(\theta,u)$ given by \eqref{eq:temps} and the uniqueness
of smooth solutions, it suffices to state that 
under the assumptions of the lemma, we have
$$\sup_{0\leq t<T}\bigl(\|\theta(t)\|_{H^{s-1}}+\|\omega(t)\|_{H^{s-1}}\bigr)<\infty.$$
First, as $\theta$ is 
transported by  the vector-field $u$
(which is lipschitz for  $s>2$ implies  $H^{s-1}\hookrightarrow L^\infty$),
 we get the following control:
$$
\|\theta(t)\|_{L^\infty}= \|\theta_0\|_{L^\infty}\quad\hbox{for
all}\quad t\in[0,T).
$$
In consequence, Inequality  \eqref{eq:apriori} ensures that 
\begin{equation}\label{eq:hs-1}
\|(\theta,\omega)(t)\|_{H^{s-1}}^2
+\nu\int_0^t\|\d_1\omega\|_{H^{s-1}}^2\,d\tau
\leq \|(\theta_0,\omega_0)\|_{H^{s-1}}^2\,
e^{(\nu^{-1}+C\|\theta_0\|_{L^\infty})t}
\,e^{C\int_0^t\|\nabla u\|_{L^\infty}\,d\tau}.
\end{equation}
On the other hand, in the appendix, it is shown that 
\begin{equation}\label{eq:log}
\|\nabla u\|_{L^\infty}\leq C\bigl(1+\|\nabla u\|_{LL}
\log(e+\|\omega\|_{H^{s-1}})\bigr).
\end{equation}
Putting together \eqref{eq:hs-1} and \eqref{eq:log}, 
we deduce that for all $t\in[0,T),$
$$\displaylines{
\log\bigl(e+\|(\theta,\omega)(t)\|_{H^{s-1}}^2\bigr)
\leq \log\bigl(e+\|(\theta_0,\omega_0)\|_{H^{s-1}}^2\bigr)
+C\bigl(1+\nu^{-1}+\|\theta_0\|_{L^\infty}\bigr)t
\hfill\cr\hfill+C\int_0^t\|\nabla u\|_{LL}
\log\bigl(e+\|(\theta,\omega)(t)\|_{H^{s-1}}^2\bigr)\,d\tau,}
$$
whence, according to  Gronwall Lemma, 
$$\log\bigl(e+\|(\theta,\omega)(t)\|_{H^{s-1}}^2\bigr)
\leq \Bigl(\log\bigl(e+\|(\theta_0,\omega_0)\|_{H^{s-1}}^2\bigr)
+C\bigl(1+\nu^{-1}+\|\theta_0\|_{L^\infty}\bigr)t\Bigr)
e^{C\int_0^t\|\nabla u\|_{LL}\,d\tau}.
$$
As the argument of $\exp$ is, by assumption, bounded for $t\in[0,T),$ 
we gather that $(\theta,\omega)\in L^\infty([0,T);H^{s-1}),$
which completes the proof of the lemma.
\end{p}
The next step involves showing that the norm used
in the previous lemma is controlled by the system. 
In fact, we shall state a slightly more accurate result:
\begin{lemma}\label{l:propagation}
Let $(\theta,u)$ be a solution to $\eqref{eq:b-vitesse}$
in $\cC([0,T);H^{2}\times H^3).$ 
There exists a continuous function $f:\R_+\rightarrow\R_+$  
depending only (continuously) on  $\nu,$ $\|u_0\|_{L^2},$
$\|\theta_0\|_{L^2\cap L^\infty}$
and $\|\omega_0\|_{\sqrt L}$
such that 
$$
\int_0^t\|\nabla u\|_{\sqrt L}^2\,dt\leq f(t)\quad\hbox{for all }\ t\in[0,T).
$$
\end{lemma}
\begin{p}
Let us first notice that, because $\div u=0$ and $u$ is Lipschitz, we have
\begin{equation}\label{eq:thetaLp}\|\theta(t)\|_{L^p}=
\|\theta_0\|_{L^p}
\quad\hbox{for all}\quad p\in[2,\infty]\quad\hbox{and}\quad
t\in[0,T).
\end{equation}
In order to get a  control 
of $\nabla u$ in $L^2([0,T);\sqrt L),$
we are going to state that 
\begin{equation}\label{eq:omegaLp}
\|\omega(t)\|_{\sqrt L}^2\leq \|\omega_0\|_{\sqrt L}^2+\frac{t}{2\nu}
\|\theta_0\|_{L^2\cap L^\infty}^2.
\end{equation}
For showing that, one may multiply the vorticity equation with
 $|\omega|^{p-2}\omega$ and perform a space integration.
As our  hypotheses on the solution entail that
$\omega\in\cC^1([0,T);H^1)$ and thus  $\omega\in\cC^1([0,T);L^p)$ for
all $p\in[2,+\infty),$ we get
$$
\begin{array}{lll}
\Frac{1}{p}\Frac{d}{dt}\Int|\omega|^p\,dx
+(p\!-\!1)\nu\!\Int\!|\partial_1\omega|^2|\omega|^{p-2}\,dx
&\leq (p\!-\!1)\Int|\theta||\partial_1\omega||\omega|^{p-2}\,dx,\\[1ex]
&\leq
(p\!-\!1)\nu\!
\Int|\partial_1\omega|^2|\omega|^{p-2}\,dx
+\frac{p\!-\!1}{4\nu}\!\int\!|\theta|^2|\omega|^{p-2}\,dx\\[1ex]&\leq
(p\!-\!1)\nu\!\Int\!|\partial_1\omega|^2|\omega|^{p-2}\,dx
+\frac{p\!-\!1}{4\nu}\,\|\theta\|_{L^p}^2\,
\|\omega\|_{L^p}^{p\!-\!2}.
\end{array}
$$
Therefore,
$$
\frac{d}{dt}\|\omega\|_{L^p}^2\leq
\biggl(\frac{p-1}{2\nu}\biggr)\|\theta\|_{L^p}^2
$$ 
hence, by virtue of  \eqref{eq:thetaLp},
$$\|\omega(t)\|_{L^p}^2\leq
\|\omega_0\|_{L^p}^2+\biggl(\frac{p-1}{2\nu}\biggr)\|\theta_0\|_{L^p}^2\,t,
$$
whence \eqref{eq:omegaLp}.
\smallbreak
This \emph{does not} imply that $\nabla u\in L^\infty([0,T);\sqrt L)$ 
for the classical result on Calderon-Zygmund operators (see e.g. \cite{C}, Chap. 3) 
gives only that 
\begin{equation}\label{eq:bs}
\|\nabla u(t)\|_{L^p}\leq Cp\|\omega(t)\|_{L^p}\quad\hbox{for all}\quad p\in[2,\infty).
\end{equation}
However, because 
 $\partial_{1}\omega=\Delta u^2$  and
  $\partial_{1}u^1=-\partial_{2}u^2,$  Inequalities 
\eqref{eq:energie} and \eqref{eq:omegaL2}
entail that
\begin{equation}\label{eq:diuj}
\sqrt\nu\|\d_iu^j\|_{L_t^2(H^1)}\leq \|u_0\|_{H^1}+\biggl(t+\sqrt{\frac t\nu}\biggr)
\|\theta_0\|_{L^2}\quad\hbox{for }\ t\in[0,T)\ \hbox{ and }\ (i,j)\not=(2,1).
\end{equation}
By virtue of Lemma \ref{l:injection} (see the appendix), we thus get 
the desired bound for all the components of $\nabla u$ except $\d_2u^1.$
In order to get a suitable bound for  $\partial_2 u^1,$ one may use
the fact that
 $\partial_2 u^1=\partial_1 u^2-\omega.$
 Putting together Inequalities  \eqref{eq:omegaLp}
and \eqref{eq:diuj}, it is now easy to conclude
 \end{p}
 \smallbreak\noindent{\it Proof of Theorem \ref{th:resultat3}~:}
 For the sake of simplicity, we restrict ourselves to the case $s\geq3$ 
 so that one may use Lemma \ref{l:propagation}. 
 The case  $2<s<3$ easily follows from the case $s\geq3$: it is only a 
 matter of smoothing out the initial data then pass to the limit. 
 \smallbreak
  So let us assume from now on that $s\geq3$ and let us
  denote by $(\theta,u)$  the maximal solution supplied by Proposition \ref{p:local}, and by
   $T^*$ the lifespan of $(\theta,u).$ If we assume (by contradiction) that $T^*$ is finite
   then Proposition \ref{l:propagation} ensures that 
 \begin{equation}\label{eq:borne1}
 \nabla u\in L^2([0,T^*);\sqrt L).
 \end{equation}
 Remark that the space $\sqrt L$ is continuously embedded in the space 
    $LL^{\frac12}$  defined by 
   \begin{equation}\label{eq:LL1/2}
LL^{\frac12}:=\biggl\{f\in {\mathcal S}'\,/\,
\|f\|_{LL^{\frac12}}:=\sup_{q\geq0}
\frac{\|S_qf\|_{L^\infty}}{\sqrt{q+1}}<\infty\biggr\}.
\end{equation}
 Indeed, thanks to  Bernstein inequality, there exists a constant
 $C$ such that for all $N\in\N$ and  $p\in[2,\infty[,$ we have
\begin{equation}\label{eq:injection1}
\|S_N\nabla u\|_{L^\infty}\leq C2^{\frac{2N}{p}}\|\nabla u\|_{L^p}
\leq C2^{\frac{2N}{p}}\sqrt {p-1}\,\|\nabla u\|_{\sqrt L}.\end{equation}
 If we choose  $p=N+2$ then we get
\begin{equation}\label{eq:injection}
\|S_N\nabla u\|_{L^\infty}\leq C \sqrt{N+1}\,\|\nabla u\|_{\sqrt L},
\end{equation}
whence the desired embedding.
\smallbreak
It is obvious that  $LL^{\frac12}\hookrightarrow LL.$
Resorting to \eqref{eq:borne1}, we thus get $\nabla u\in L^1([0,T^*);LL)$
and Lemma \ref{l:explosion} ensures that the solution $(\theta,u)$ may
be continued beyond  $T^*.$
This contradicts the definition of $T^*.$
 \hfill\rule{2.1mm}{2.1mm}


\subsection{Global well-posedness for rough data}

 In this  section, we want to state  global  existence with 
 uniqueness for  a class of data as large as possible.
 Having in mind the previous subsection, it seems reasonable to require that 
$\theta_0\in L^2\cap L^\infty,$
 that   $u_0\in H^1$  and that 
$\omega_0\in\sqrt L.$
As those regularity assumptions are (formally) conserved
by the system during the evolution, 
we thus expect to get a global solution 
 $(\theta,u)$  such that  $\nabla u\in L^2_{loc}(\R_+;\sqrt L).$
This would imply that $u\in L^2_{loc}(\R_+;{\rm LogLip}^{\frac12})$ 
where  ${\rm LogLip}^{\frac12}$ stands for the set of bounded functions $f$ such that 
$$
\sup_{\substack{x\not=y\\|x-y|\leq1/2}}
\frac{|f(y)-f(x)|}{|x-y|\log^{\frac12}(|x-y|^{-1})}<\infty.
$$
The above inequality is an obvious corollary of \eqref{eq:injection}
and of Proposition 2.107 in \cite{BCD}. 

Even though  the vector-field $u$ fails to be lipschitz, 
it  has enough regularity so that
 we have $\theta\in\cC(\R_+;H^{s-\e})$ for all  $\epsilon>0$
  if we start from  $\theta_0$in $H^s$ 
for some $s\in]-1,1]$ 
(this is in fact a consequence of Theorem 3.12 in  \cite{DANCHI}). 
These plain observations will lead us to the following
statement which obviously contains Theorem \ref{th:resultat1}:
\begin{theorem}\label{th:globalfaibleregulier}
Let  $\theta_0\in L^2\cap L^\infty,$ and   $u_0\in H^1$ with $\div u_0=0$
and  $\omega_0\in\sqrt L.$
 Then System~$\eqref{eq:b-vitesse}$ admits a global solution 
 $(\theta,u)$ such that 
\begin{equation}\label{eq:regu}
\begin{array}{c}\theta\in\cC_b(\R_+;L^2)\cap \cC_w(\R_+;L^\infty)\cap L^\infty(\R_+;L^\infty),
\\[1.5ex]
u\in\cC_w(\R_+;H^1),\quad\! u^2\in L^2_{loc}(\R_+;H^2),
\!\quad
\omega\in L^\infty_{loc}(\R_+;\sqrt L),\quad\! 
\nabla u\!\in\! L^2_{loc}(\R_+;\sqrt L).\end{array}
\end{equation}
If in addition  $\theta_0\in H^s$ for some real number
$s\in(0,1]$ then
$\theta\in\cC(\R_+;H^{s-\e})$ for all $\e>0.$

Finally, if $s>1/2$ then the  solution is unique.
\end{theorem}
\begin{p} 
The uniqueness is a consequence 
of Proposition
\ref{p:unicite} below so let us focus on the
existence part of the statement.  
To achieve it, one may smooth out the initial data
$(\theta_0,u_0)_{n\in\N}$ so as to get a sequence 
 $(\theta_{0,n},u_{0,n})_{n\in\N}$ of $H^\infty$ 
 functions which tends to $(\theta_0,u_0)$ in 
 (say) $L^2\times H^1.$
 Resorting to Theorem \ref{th:resultat3}, we get a sequence
 $(\theta_n,u_n)_{n\in\N}$ of smooth global solutions.
Moreover, by virtue of Inequalities \eqref{eq:energie}, \eqref{eq:omegaL2}, 
\eqref{eq:thetaLp} and \eqref{eq:omegaLp}, we have
\begin{itemize}
\item $(\theta_n)_{n\in\N}$ bounded in $L^\infty(\R_+;L^2\cap L^\infty),$
\item $(u_n)_{n\in\N}$ bounded in $L^\infty_{loc}(\R_+;H^1),$
\item  $(u_n^2)_{n\in\N}$ bounded in $L^2_{loc}(\R_+;H^2),$
\item  $(\omega_n^2)_{n\in\N}$ bounded in  $L^\infty_{loc}(\R_+;\sqrt L).$
\end{itemize}
As explained in the proof of Lemma 
\eqref{l:propagation}, these proprieties
imply that  $(\nabla u_n)_{n\in\N}$ is bounded
in $L^2_{loc}(\R_+;\sqrt L).$

Now, taking advantage of the losing estimates
proved in  Proposition \ref{p:transport}, 
we deduce that if, in addition, $\theta_0\in H^s$ for some $s\in(0,1)$ 
then, for all $\e>0,$ 
$(\theta_n)_{n\in\N}$ is bounded in  $L^\infty_{loc}(\R_+;H^{s-\e}).$
In order to conclude to the existence part of 
the statement, one may use again 
a compactness argument  \emph{\`a la Aubin-Lions} 
as in the proof of Theorem \ref{th:globalfaible}.
\end{p}


\subsection{Uniqueness for rough data }\label{ss:uniqueness}

The difference $(\dt,\du,\dPi)$
between two solutions $(\theta_1,u_1,\Pi_1)$ and $(\theta_2,u_2,\Pi_2)$ satisfies:
\begin{equation}
\label{eq:unicite}
\begin{cases}
\partial_t\dt+\div(u_2\,\dt)=-\div(\du\, \theta_1),
\\ \partial_t\du+\div(u_2\otimes\du)+\div(\du\otimes u_1)-\nu\d_1^2\du
+\nabla\dPi=\dt\,e_2.
\end{cases}
\end{equation}
First of all, let us notice that if $\theta_1\in L^\infty([0,T];H^\gamma)$
then the right-hand side of the first equation is at most in  
$L^\infty([0,T];H^{\gamma-1}).$
As, under the hypotheses of Theorem
\ref{th:resultat3},   $u_2$  fails to be  lipschitz but has  gradient in
$L^2_{loc}(\R_+;\sqrt L),$ Proposition \ref{p:transport} will enable us to bound
   $\dt$ in   $H^{\gamma-1-\e}$
 provided that  $\div(\du\,\theta_1)$ be   bounded in $L^1([0,T];H^{\gamma-1}),$
 a condition which requires   a control over the  $L^1([0,T];L^\infty)$ norm of
  $\du.$ If the velocity equation had a full  Laplace operator
  then the resulting smoothing effect would be  strong enough 
  so as to provide us with a control over this norm. 
  We shall see that in the framework of anisotropic viscosity, 
  one can still get an appropriate bound for  $\du$
in  $L^1([0,T];L^\infty)$  provided  $\gamma>1/2.$
In short, we expect to be able to control 
$\dt$ in $L^\infty([0,T];H^{\beta-1})$ for some  $\beta\in]1/2,\gamma[,$
 $\du$ in  $L^\infty([0,T];H^\alpha)$ 
and $\d_1\du$ in  $L^2([0,T[;H^\alpha)$ 
for some  $\alpha\in]1/2,\beta[.$
This motivates the following statement 
which implies the uniqueness part of Theorem \ref{th:globalfaibleregulier}:
\begin{proposition}\label{p:unicite}
Let $(\theta_1,u_1)$ and  $(\theta_2,u_2)$ be two 
solutions of  \eqref{eq:b-vitesse} on $[0,T]\times\R^2$
with the same initial data. 
Assume that $\nabla u_i\in L^1([0,T];LL^{\frac12})$ 
and that there exists some  $\gamma\in]1/2,1[$ such that 
$$
\theta_i\in L^\infty([0,T];L^\infty\cap H^\gamma)\ \hbox{ and }\ 
u_i\in L^\infty([0,T];H^\gamma)\quad\hbox{for }\ i=1,2.
$$
Then the two solutions coincide.
\end{proposition}
\begin{p}
According to the above  heuristics, we have 
to bound  
$$\dt\ \hbox{ in }\ L^\infty([0,T];H^{\beta-1}),\quad
\du\ \hbox{ in }\  L^\infty([0,T];H^\alpha)\ \hbox{ and }\ 
\d_1\du\ \hbox{ in }\  L^2([0,T];H^\alpha)
$$
for some fixed  $(\alpha,\beta)$ such that   $1/2<\alpha<\beta<\gamma.$
\smallbreak
In order to bound $\dt,$ one may use Proposition \ref{p:transport}
with the vector-field $u_2.$
Because  $\dt(0)=0,$ we have
$$
\Vert\dt(t)\Vert_{H^{\beta-1}}
\leq C\int_0^t\Vert \div(\du\,\theta_1)\|_{H^{\gamma-1}}\quad\hbox{for all}\quad t\in[0,T]. $$
In order to bound the right-hand side, one may resort
to the following  \emph{Bony's decomposition} \cite{B}: 
\begin{equation}\label{eq:bony}
\div(\du\,\theta_1)=
\div\Bigl(T_{\du}\theta_1+R(\du,\theta_1)\Bigr)
+\sum_{i=1}^2T_{\d_i\theta_1}\du^i,
\end{equation}
where the paraproduct operator  $T$ (resp. reminder operator $R$)
is defined by 
$$
T_{f}g:=\sum_{q} S_{q-1}f\dq g\qquad
\biggl(\hbox{resp. }\quad
R(f,g):=\sum_q \dq f\bigl(\Delta_{q-1}g\!+\!\dq g\!+\!\Delta_{q-1}g\bigr)\biggr).
$$
Let us stress that the condition $\div\du=0$ has
been used in order to have  the derivative act on the left for the first two terms
of  \eqref{eq:bony}.
\smallbreak
From standard continuity results for operators $T$ and $R$
(see e.g. \cite{BCD})  we have
$$
\Vert T_{\du}\theta_1+R(\du,\theta_1)\Vert_{H^{\gamma}}
\leq C\Vert\du\Vert_{L^\infty}\Vert\theta_1\Vert_{H^{\gamma}}.$$
As for the last term, given that $\gamma-1<0,$ one can write
$$
\|T_{\d_i\theta_1}\du^i\|_{H^{\gamma-1}}
\leq C\|\nabla\theta_1\|_{H^{\gamma-1}}\|\du\|_{L^\infty}.
$$
We eventually get
\begin{equation}\label{eq:unicite1}
\|\delta\!\theta\|_{L_t^\infty(H^{\beta-1})}
\leq C \Vert\theta_1\Vert_{L_t^\infty(H^\gamma)}
\|\du\|_{L_t^1(L^\infty)}.
\end{equation}
In order to bound  $\du,$ one may use Proposition 
\ref{p:stokespartiel}.
We get for all  $t\in[0,T],$
$$
\|\du\|_{L^\infty_t(H^\alpha)}+
\|\d_1\du\|_{L_t^2(H^\alpha)}
\leq C\bigl(\|\dt\|_{L^2_t(H^{\beta-1})}
+\|\du\cdot\nabla u_1\|_{L^2_t(H^{\beta-1})}\bigr)
$$
for some constant $C$ depending only on 
 $\alpha,$ $\beta,$ $\nu$ and  $u_2.$
 \smallskip
Using again the Bony decomposition
and arguing exactly as for proving  \eqref{eq:unicite1},
we get
$$
\|\du\cdot\nabla u_1\|_{H^{\beta-1}}
\leq C\|\du\|_{L^\infty}\|u_1\|_{H^\beta}.
$$
Therefore, given that $u_1\in L^\infty([0,T];H^\beta),$ 
\begin{equation}\label{eq:unicite2}
\|\du\|_{L^\infty_t(H^\alpha)}+
\|\d_1\du\|_{L_t^2(H^\alpha)}
\leq C\bigl(\|\dt\|_{L^2_t(H^{\beta-1})}
+\|\du\|_{L^2_t(L^\infty)}\bigr).
\end{equation}
In order to complete the proof of the proposition,
it is only a matter of showing that 
$\|\du\|_{L^1_t(L^\infty)}$
may be bounded in terms of 
$\|\du\|_{L^\infty_t(H^\alpha)}$ and of  
$\|\d_1\du\|_{L_t^2(H^\alpha)}.$
This is the only point where the assumption $\alpha>1/2$ (and thus $\gamma>1/2$)
is going to play a role.
First of all, thanks to the trace theorem, one may write (with obvious notation)
 $$
H^\alpha(\R^2)\hookrightarrow L^\infty(\R_{x_2};H^{\alpha-\frac12}(\R_{x_1})).
$$
Therefore, one may write 
\begin{equation}\label{eq:unicite3}
\Vert\du\Vert_{L^\infty_{x_2}(H^{\alpha-\frac12}_{x_1})}
\leq C\|\du\|_{H^\alpha}\quad\hbox{and}\quad
\Vert\d_1\du\Vert_{L^\infty_{x_2}(H^{\alpha-\frac12}_{x_1})}
\leq C\|\d_1\du\|_{H^\alpha}.
\end{equation}
As for all  $\alpha\in]0,1[,$ 
Gagliardo-Nirenberg  inequality implies that 
$$
\|\du(\cdot,x_2)\|_{L^\infty(\R)}
\leq C\|\du(\cdot,x_2)\|_{H^{\alpha-\frac12}(\R)}^\alpha
\|\d_1\du(\cdot,x_2)\|_{H^{\alpha-\frac12}(\R)}^{1-\alpha}\quad\hbox{for all }\ x_2\in\R,
$$
we have, by combination with \eqref{eq:unicite3}, 
$$
\|\du\|_{L^\infty(\R^2)}
\leq C\|\du\|_{H^\alpha(\R^2)}^\alpha\|\d_1\du\|_{H^\alpha(\R^2)}^{1-\alpha}.
$$
Coming back to \eqref{eq:unicite1} and  \eqref{eq:unicite2},
we deduce that for some constant $C$ depending only on 
$\nu,$ $T$  and on the  norms of  $(\theta_1,u_1)$ and  $(\theta_2,u_2),$ we have
$$\displaylines{
\|\delta\!\theta\|_{L_t^\infty(H^{\beta-1})}
\leq Ct^{\frac12+\frac\alpha2}\dU(t)\cr
\dU(t)\leq C\Bigl(t^{\frac12}
\|\delta\!\theta\|_{L_t^\infty(H^{\beta-1})}
+t^{\frac\alpha2}\dU(t)\Bigr)}
$$
with 
$$
\dU(t):=\|\du\|_{L_t^\infty(H^\alpha)}
+\|\d_1\du\|_{L_t^2(H^\alpha)}.
$$
Inserting the first inequality in the second one, 
one may conclude that  $\du\equiv0$ (and thus $\dt\equiv0$) 
on a suitably small time interval.
Finally, let us notice that our assumptions on the solutions ensure
that 
$\dt\in\cC([0,T];H^{\beta-1})$ and  $\du\in\cC([0,T];H^\alpha).$
Using a classical connectivity  argument, it is now easy to get
the uniqueness on the whole interval  $[0,T].$\end{p}


\section{The case of an horizontal diffusivity} \label{s:b-temperature}

This section is devoted to the study of System \eqref{eq:b-temperature}.
In other words,  in contrast with the previous section, 
we now assume that the velocity satisfies 
the incompressible Euler equation with buoyancy force
whereas  the temperature experiences diffusion \emph{in the horizontal
variable only}. 

We aim at stating various global existence results
for arbitrarily large data. 
More precisely, we first prove that any data  
$\theta_0\in L^2$ and $u_0\in H^1$ with $\div u_0=0$
generates a global weak solution with finite energy. 
The rest of this section is mainly devoted to the proof
of Theorem \ref{th:resultat2}.
As a first step, in subsection
 \ref{ss:H1}, we state  $H^1$ a priori estimates for the
 la temperature. In the next subsection, we prove 
a uniqueness result for a large class of solutions.
 As this uniqueness result requires in particular that 
$\nabla\theta\in L^1_{loc}(\R_+;L^\infty)$ and that 
$\nabla u\in L^1_{loc}([0,T];L),$ our next
task amounts to finding additional regularity 
conditions on the data which may be propagated globally by the
system. 
It turns out that it is possible to propagate
some  \emph{anisotropic} Sobolev regularity over the temperature,
and thus to complete the proof of Theorem~\ref{th:resultat2}.
 
 
\subsection{Global weak solutions: the  case
 $\theta_0\in L^2$ and  $u_0\in H^1$}

Let us first derive the formal  energy estimates
for System \eqref{eq:b-temperature}
in the case   $\theta_0\in L^2$ and $u_0\in H^1.$
First, multiply $\eqref{eq:b-temperature}_1$ 
by  $\theta$ and integrate over $[0,t]\times\R^2.$  
We get
\begin{equation}\label{eq:theta1}
\|\theta (t)\|_{L^2}^2+2\kappa\int_0^t\|\partial_1\theta(s)\|_{L^2}^2\,ds
\leq \|\theta_0\|^2_{L^2}.\end{equation}
Combining this with  the standard energy estimate for $u$ yields
\begin{equation}\label{eq:u1}
\|u(t)\|_{L^2}\leq\|u_0\|_{L^2}+t\|\theta_0\|_{L^2}.
\end{equation}
In order to get a  $H^1$ bound for the velocity, one may 
consider the vorticity equation: 
$$\partial_t\omega+u\cdot\nabla\omega=\partial_1\theta.$$
Multiplying by $\omega$ then integrating with respect to 
the space variable, we find that
$$\frac 12\frac{d}{dt}\|\omega\|^2_{L^2}\leq 
\|\omega\|_{L^2}\|\partial_1\theta\|_{L^2}$$
whence,  
\begin{eqnarray}
\|\omega(t)\|_{L^2}&\leq& \|\omega_0\|_{L^2}
+\int_0^t\|\partial_1\theta\|_{L^2}\,ds,\nonumber\\\label{eq:omega1}
&\leq&
\|\omega_0\|_{L^2}+\sqrt{\frac t{2\kappa}}\|\theta_0\|_{L^2}.
\end{eqnarray} 
Now, using a Friedrichs method quite similar 
to that of the proof of Theorem \ref{th:globalfaible}, 
we easily get the following statement:
\begin{theorem}\label{th:resultat4} 
Let $\theta_0\in L^2$ and  $u_0\in H^1$ with $\div u_0=0.$
Then System $\eqref{eq:b-temperature}$ with data 
$(\theta_0,u_0)$
has a global solution  $(\theta,u)$  such that
$$
\theta\in\cC_w(\R_+;L^2),
\quad \partial_1\theta\in L^2(\R_+;L^2)
\quad\!\hbox{and}\!\quad
u\in\cC(\R_+;H^1).
$$
\end{theorem}


\subsection{$H^1$ a priori estimates for the temperature}\label{ss:H1}
In the present  paragraph, we show that one may get (at least formally)
a  global control over both the  $H^1$ norm of  $\theta$  and of  $u.$
\smallbreak
To start with, let us point out that Inequalities \eqref{eq:theta1},
\eqref{eq:u1} and \eqref{eq:omega1}  provide us with a bound 
for $\theta$ in
$L^\infty(\R_+;L^2),$ 
for $\d_1\theta$ in 
$L^2(\R_+;L^2)$ and for $u$ in  $L^\infty_{loc}(\R_+;H^1).$ 
We claim that if we assume in addition that 
$\nabla\theta_0\in L^2$ then one may bound $\nabla\theta$ in
$L^\infty_{loc}(\R_+;L^2).$ 
Indeed, applying operator $\d_i$ ($i=1,2$) to the equation
satisfied
by  $\theta$ yields  
$$\partial_t\partial_i\theta+u\cdot\nabla\partial_i\theta+\partial_iu\cdot\nabla
\theta-\kappa\partial_1^2\d_i\theta=0.$$ 
Let us multiply this equality by $\d_i\theta,$
integrate over $\R^2$ then add up the  equalities
for $i=1,2.$
Integrating by parts where needed and using the fact that $\div u=0,$
we easily find that 
  \begin{equation}\label{eq:nablatheta}
\Frac
12\Frac{d}{dt}\|\nabla\theta\|^2_{L^2}
+\kappa\|\partial_1\nabla\theta\|^2_{L^2}+\sum_{1\leq i,j\leq2}
\Int \partial_i\theta\,\d_ju^i\,\d_j\theta\,dx=0.
\end{equation}
For $(i,j)\not=(2,2),$ the terms in the above summation 
are easy to handle. 
Indeed, taking advantage of the anisotropic H\"older inequality, 
one can write
$$
\biggl|\Int \partial_i\theta\,\d_ju^i\,\d_j\theta\,dx\biggr|\leq
\|\nabla u\|_{L^2}\|\d_1\theta\|_{L^2_{x_1}(L^\infty_{x_2})}
\|\nabla\theta\|_{L^\infty_{x_1}(L^2_{x_2})}.
$$
Let us admit the following two inequalities (the proof
of which is postponed in the appendix): 
\begin{equation}\label{eq:anisobound}
\|f\|_{L^2_{x_1}(L^\infty_{x_2})}\leq 
C\|f\|_{L^2}^{\frac 12}\|\partial_2f\|_{L^2}^{\frac 12}
\quad\hbox{and}\quad\|f\|_{L^\infty_{x_1}(L^2_{x_2})}\leq 
C\|f\|^{\frac 12}_{L^2}\|\partial_1f\|^{\frac 12}_{L^2}.\end{equation}
Applying these inequalities to  $\d_1\theta$ and to  $\nabla\theta,$ 
and using the fact that  $\|\nabla u\|_{L^2}=\|\omega\|_{L^2},$
we deduce that
$$
\biggl|\Int \partial_i\theta\,\d_ju^i\,\d_j\theta\,dx\biggr|\leq
C\|\omega\|_{L^2}\|\nabla\theta\|_{L^2}\|\d_1\nabla\theta\|_{L^2}
\quad\hbox{if}\quad (i,j)\not=(2,2).
$$
In order to bound the term corresponding to $(i,j)=(2,2),$ one may 
use the fact that $\d_2u^2=-\d_1u^1$ and integrate by parts. We get
$$
\int\d_2u^2 (\d_2\theta)^2\,dx=-\int\d_1u^1(\d_2\theta)^2\,dx
=2\int u^1\d_2\theta\,\d_1\d_2\theta\,dx.
$$
Therefore, thanks to the anisotropic H\"older inequalities and
to \eqref{eq:anisobound}, 
$$\begin{array}{lll}
\biggl|\Int\d_2u^2 (\d_2\theta)^2\,dx\biggr|&\leq&
2\|\d_1\d_2\theta\|_{L^2}\|u^1\|_{L^2_{x_1}(L^\infty_{x_2})}
\|\d_2\theta\|_{L^\infty_{x_1}(L^2_{x_2})},\\[1ex]
&\leq&C\|\d_1\d_2\theta\|_{L^2}\|u\|_{L^2}^{\frac12}\|\omega\|_{L^2}^{\frac12}
\|\d_2\theta\|_{L^2}^{\frac12}\|\d_1\d_2\theta\|_{L^2}^{\frac12},\\[1ex]
&\leq&C\|u\|_{L^2}^{\frac12}\|\omega\|_{L^2}^{\frac12}
\|\d_2\theta\|_{L^2}^{\frac12}\|\d_1\d_2\theta\|_{L^2}^{\frac32}.
\end{array}
$$
So finally,  Young inequality leads to  
$$
\biggl|\sum_{1\leq i,j\leq2}
\Int \partial_i\theta\,\d_ju^i\,\d_j\theta\,dx\biggr|
\leq \frac\kappa2\|\d_1\nabla\theta\|_{L^2}^2
+\frac C\kappa\|\omega\|_{L^2}^2\biggl(1+\frac{\|u\|_{L^2}^2}{\kappa^2}\biggr)
\|\nabla\theta\|_{L^2}^2.
$$
Plugging this inequality in \eqref{eq:nablatheta}
and using  Gronwall lemma, we end up with

$$\|\nabla\theta(t)\|_{L^2}^2+\kappa
\int_0^t\|\partial_1\nabla\theta(s)\|^2_{L^2}\,ds
\leq \|\nabla\theta_0\|_{L^2}^2\exp\biggl\{\frac C\kappa\int_0^t
\|\omega\|_{L^2}^2
\biggl(1+\frac{\|u\|_{L^2}^2}{\kappa^2}\biggr)\,d\tau\biggr\}.
$$
Putting together \eqref{eq:theta1}, \eqref{eq:u1} and \eqref{eq:omega1},
we conclude that 
\begin{equation}
\label{eq:estimation:h1}
\|\theta(t)\|_{H^1}^2+\kappa
\int_0^t\|\partial_1\theta(s)\|^2_{H^1}\,ds
\leq C(t,\kappa,\theta_0,u_0)
\end{equation}
with 
$\,\displaystyle{C(t,\kappa,\theta_0,u_0):=\|\theta_0\|_{H^1}^2\exp
\bigg\{\frac{Ct}\kappa\biggl(\|\omega_0\|_{L^2}^2+\frac t\kappa
\|\theta_0\|_{L^2}^2\biggr)\biggl(1+\frac{\|u_0\|_{L^2}^2
+t^2\|\theta_0\|_{L^2}^2}{\kappa^2}\biggr)\bigg\}\cdotp}$


\subsection{A uniqueness result}

In this section, we establish a uniqueness result
for System \eqref{eq:b-temperature} under ``minimal" 
assumptions. In order to motivate those assumptions, let us remind 
that in the isotropic case (that is with a full Laplacian 
in the temperature equation) which has been
investigated in  \cite{DP3}, uniqueness
is true in the class of $C^{0,1}(\R_+;L^2)$ solutions which satisfy in 
addition 
\begin{equation}\label{eq:iso}
\nabla\theta\in L^1_{loc}(\R_+;L^\infty)\quad\hbox{and}\quad
\nabla u\in L^1_{loc}(\R_+;L).
\end{equation}
As in the case that we now consider the smoothing effect over the temperature 
is obviously \emph{weaker}, 
we expect the conditions leading to uniqueness to be \emph{stronger}
than \eqref{eq:iso}. 
We shall prove the following result:
\begin{proposition}\label{p:uniqueness}
Let $(\theta_1,u_1)$ and  $(\theta_2,u_2)$ be two solutions of  
$\eqref{eq:b-temperature}$ with the same data. Assume 
that both solutions belong to 
$L^\infty([0,T];H^1)\cap C^{0,1}([0,T];L^2)$ and that, in addition,
$\d_1\theta_2\in L^2([0,T];H^1)$ and  $\nabla u_2\in L^1([0,T];L).$ 
Then 
$(\theta_1,u_1)\equiv(\theta_2,u_2)$ on $[0,T]\times\R^2.$ \end{proposition}
\begin{p}  
With the usual  notation,   
 $(\dt, \du)$ satisfies: 
$$
\begin{cases}
\partial_t\dt+u_1\cdot\nabla\dt+\div(\du\theta_2)
-\kappa\partial_1^2\dt=0\\
\partial_t \du+u_1\cdot\nabla \du+\du\cdot\nabla u_2=-\nabla
\dPi+\dt\,e_2. \end{cases}
$$
{}From a standard energy method, we get
\begin{eqnarray}\label{eq:dtheta}&\Frac
12\frac{d}{dt}\|\dt\|_{L^2}^2+\kappa\|\partial_1\dt\|^2_{L^2}\leq
\bigg|\int\div(\theta_2\du)\dt\,dx\bigg|,\\\label{eq:dvitesse}&\Frac
12\frac{d}{dt}\|\du\|_{L^2}^2\leq \bigg|\int\du\cdot\nabla u_2\:\dt\,dx
\bigg|+\bigg|\int\dt \du^2\,dx\bigg|.\end{eqnarray} 
In order to bound the right-hand side of \eqref{eq:dtheta}, one may write
  \begin{equation}\label{eq:dtheta0}\int \div(\theta_2 \du)\dt=-\int
\theta_2 \du^1\partial_1\delta \theta\,dx-\int \theta_2
\du^2\partial_2\dt\,dx.\end{equation} 
The first term is easy to deal with:
using  Cauchy-Schwarz inequality, we get 
\begin{equation}\label{eq:dtheta1}
\biggl|\Int\theta_2\du_1\partial_1\dt\,dx\biggr|\leq
\|\du\|_{L^2}\|\theta_2\|_{L^\infty}\|\partial_1\dt\|_{L^2}. \end{equation}
Next, applying the following inequality  (see the proof in the appendix)
 \begin{equation}\label{eq:anisobound1}
\|\theta_2\|_{L^\infty}\leq C\|\theta_2\|_{L^2}^{1/4}\|\partial_1\theta_2\|_{L^2}^{1/4}
\|\partial_2\theta_2\|_{L^2}^{1/4} 
\|\partial_1\partial_2\theta_2\|_{L^2}^{1/4},
\end{equation}
and using Young inequality, we find that
\begin{equation}\label{eq:dtheta2}\left|\Int\theta_2\du_1\partial_1\dt\,dx\right|\leq
\frac{C}{\kappa}\|\theta_2\|_{H^1}^{3/2} 
\|\partial_1\partial_2\theta_2\|_{L^2}^{1/2}
\|\du\|_{L^2}^2+\frac{\kappa}{6}\|\partial_1\dt\|_{L^2}^{2}. 
\end{equation}
The second term of  \eqref{eq:dtheta0} is more intricate.
If we integrate by parts and use the fact that
$\div\du=0,$ we get 
 $$ \begin{array}{lll}
-\Int\theta_2 \delta
u^2\partial_2\dt\,dx&=&\Int\partial_2\theta_2 \du^2\: \dt\,dx+\int\theta_2\partial_2\du^2\dt\,dx,\\
&=&\Int\partial_2\theta_2\du^2\:\dt\,dx-\int\theta_2\partial_1\du^1\dt\,dx,\\
&=&A_1+A_2+A_3\end{array}$$
with
$$ A_1:=\Int\partial_2\theta_2\:\du^2\:\dt\,dx,\quad
A_2:=\int\theta_2\du^1\partial_1 \dt\,dx\ \hbox{ and }\ 
A_3:=\int\partial_1\theta_2\du^1\dt\,dx.
$$
The term $A_2$ may be bounded according to 
\eqref{eq:dtheta2}.  In order to bound $A_3,$ we use the anisotropic H\"older
inequality and  \eqref{eq:anisobound}. This leads to 
$$\begin{array}{lll}|A_3|&\leq&
\|\dt\|_{L^\infty_{x_1}(L^2_{x_2})}\|\partial_1\theta_2\|_{L^2_{x_1}(L^\infty_{x_2})} \|\du\|_{L^2},\\[2ex]&\leq&
C\|\dt\|_{L^2}^{1/2}\|\partial_1\dt\|_{L^2}^{1/2}
\|\partial_1\theta_2\|_{L^2}^{1/2}\|\partial_1\partial_2\theta_2\|_{L^2}^{1/2}
\|\du\|_{L^2},
\end{array}$$
whence, resorting again to Young inequality, 
\begin{equation}\label{eq:dtheta3}
|A_3|\leq\frac{C}{\kappa}\|\delta
\theta\|_{L^2}^2\|\partial_1\theta_2\|_{L^2}^2
\|\partial_1\partial_2\theta_2\|_{L^2}^2+\|\du\|_{L^2}^2
+\frac{\kappa}{6}\|\partial_1\dt\|_{L^2}^2.
\end{equation}
The  term $A_1$ is the most difficult to deal with.
To get an appropriate bound, let us first notice that, as  $\div\du=0,$ 
we may write 
$$\du^2=(1-\partial_2^2)^{-1}(\du^2)+(1-\partial_2^2)^{-1}\partial_2\partial_1 \du^1.$$ 
Therefore, integrating by parts, we get
$A_1=A_1^1+A_1^2+A_1^3$ with   $$\begin{array}{lll}
A_1^1&:=&\Int(1-\partial_2^2)^{-1}(\du^2)
\:\partial_2\theta_2\: \dt\,dx,\\
A_1^2&:=&-\Int \partial_2(1-\partial_2^2)^{-1}(\du^1)\:
\partial_1\partial_2\theta_2\:\dt\,dx,\\
 A_1^3&:=&-\Int\partial_2(1-\partial_2^2)^{-1}(\du^1)\:\partial_2\theta_2 \:
\partial_1\dt\,dx.
\end{array}$$
First of all, we have  
$$|A_1^3|\leq \|\partial_2(1-\d_2^2)^{-1}\du\|_{L^2_{x_1}(L^\infty_{x_2})}
\|\partial_2\theta_2\|_{L^\infty_{x_1}(L^2_{x_2})}
\|\partial_1\dt\|_{L^2}.$$ 
Taking advantage of  \eqref{eq:anisobound}, we get
$$\displaylines{\|\partial_2(1-\partial_2^2)^{-1}\du^1\|_{L^2_{x_1}(L^\infty_{x_2})}\leq
C\|\partial_2(1-\partial_2^2)^{-1}(\du^1)\|_{L^2}^{\frac
12}\|\partial_2^2(1-\partial_2^2)^{-1}\du^1\|_{L^2}^{\frac 12}\leq
C\|\du\|_{L^2}\cr
\|\partial_2\theta_2\|_{L^\infty_{x_1}(L^2_{x_2})}\leq
C\|\partial_2\theta\|_{L^2}^{\frac
12}\|\partial_1\partial_2\theta\|_{L^2}^{\frac 12}.}
$$ 
In consequence, thanks to Young inequality, we have
 $$
|A_1^3|\leq\frac C\kappa\|\partial_2\theta\|_{L^2}
\|\partial_1\partial_2\theta\|_{L^2}\|\du\|_{L^2}^2
+\frac{\kappa}{6}\|\partial_1\delta \theta\|_{L^2}^2.
$$
To deal with  $A_1^2,$ one may write that, by virtue of
 \eqref{eq:anisobound} and of Young inequality
$$\begin{array}{lll}
|A_1^2|&\leq&
\|\partial_2(1-\partial_2^2)^{-1}\du^1\|_{L^2_{x_1}(L^\infty_{x_2})}
\|\partial_1\partial_2\theta_2\|_{L^2}
\|\dt\|_{L^\infty_{x_1}(L^2_{x_2})},\\[2ex]&\leq& C\|\du\|_{L^2}
\|\partial_1\partial_2\theta_2\|_{L^2}\|\dt\|_{L^2}^{\frac
12}\|\partial_1\dt\|_{L^2}^{\frac 12},\\[2ex]&\leq& \frac C\kappa
\|\du\|_{L^2}^2\|\partial_1\partial_2\theta_2\|_{L^2}^2
+\frac{3\kappa}{2}\|\dt\|_{L^2}^2+\frac{\kappa}{6}\|\partial_1\dt\|^2_{L^2}.
\end{array}$$
Finally, for   $A_1^1$ we have
$$\begin{array}{lll}
|A_1^1|&\leq& \|(1-\d_2^2)^{-1}(\du_2)\|_{L^2_{x_1}(L^\infty_{x_2})}
\|\d_2\theta_2\|_{L^2}\|\dt\|_{L^\infty_{x_1}(L^2_{x_2})},\\[1ex]&\leq& C
\|\du\|_{L^2}\|\partial_2\theta_2\|_{L^2}\|\dt\|_{L^2}^{\frac
12}\|\partial_1\dt\|_{L^2}^{\frac 12},\\[1ex]&\leq&
\frac C\kappa\|\du\|_{L^2}^2\|\d_1\theta_2\|_{L^2}^2+\frac{3\kappa}2\|\dt\|_{L^2}^2
+\frac\kappa{6}\|\d_1\dt\|_{L^2}^2.
\end{array}$$
Putting together all the previous inequalities, we conclude that
\begin{equation}\label{eq:dtheta4}
|A_1|\leq  \frac{\kappa}{2}\|\partial_1\dt\|_{L^2}^2+3\kappa\|\dt\|_{L^2}^2+
\frac C\kappa\|\d_1\theta_2\|_{H^1}^2\|\du\|_{L^2}^2.
\end{equation}
Now, inserting Inequalities \eqref{eq:dtheta2},
\eqref{eq:dtheta3} and \eqref{eq:dtheta4} in \eqref{eq:dtheta0},
we deduce that there exists  an integrable function $f_2$ over 
$[0,T]$  depending only on $(\theta_2,u_2)$
and on $\kappa$ such that
 \begin{equation}\label{eq:dtheta6}
\frac12\frac{d}{dt}\|\dt\|_{L^2}^2\leq f_2(t)\,\|(\dt,\du)\|_{L^2}^2.
\end{equation}
Adapting the well-known Yudovich's argument  (see \cite{PG} and
\cite{Yudo}), it is now easy to complete the proof of uniqueness.
Indeed, from Inequality
\eqref{eq:dvitesse}, we get for all $p\in[2,\infty[,$
\begin{equation}\label{eq:dtheta7}\frac 12\frac{d}{dt}\|\du\|_{L^2}^2\leq
\|\nabla u_2\|_{L^p}
\|\du\|_{L^\infty}^{\frac2p}\|\du\|_{L^2}^{2-\frac2{p}}+\|\dt\|_{L^2}\|\du\|_{L^2}.
 \end{equation}  
Setting  $X_\e(t):=\sqrt{\|(\dt,\du)(t)\|_{L^2}^2+\e^2}$
for $\e>0,$ and using  \eqref{eq:dtheta6} and \eqref{eq:dtheta7}, we obtain
  $$
\frac d{dt}X_\e\leq p\|\nabla u_2\|_{L}\|\du\|_{L^\infty}^{\frac2p}
X_\e^{1-\frac2p}
+(\frac12+f_2)X_\e.
$$
Now, if we set
$Y_\e=X_\e\exp\bigl(-\int_0^t(\frac12+f_2(\tau))\,d\tau)\bigr),$ we have
$$
\frac2pY_\e^{\frac2p-1}\frac d{dt}Y_\e\leq  2\|\nabla
u_2\|_{L}\|\du\|_{L^\infty}^{\frac2p},  $$
 whence
 $$
Y_\e(t)\leq\biggl(\e^{\frac2p}
+2\int_0^t\|\nabla
u_2\|_{L}\|\du\|_{L^\infty}^{\frac2p}\,d\tau\biggr)^{\frac p2}.
$$
Having  $\e$ tend to  $0,$ we discover that for all $t\in\R_+,$
$$
\|(\dt,\du)(t)\|_{L^2}^2
\leq \|\du\|_{L_t^\infty(L^\infty)}^2 \biggl(2\int_0^t\|\nabla
u_2\|_{L}\,d\tau\biggr)^p.
$$
By Sobolev embedding and thanks to \eqref{eq:bs} with $p=4,$ we have
\begin{equation}\label{eq:uLinfty}
\|\du\|_{L^\infty}\leq C\bigl(\|\du\|_{L^2}+\|\delta\!\omega\|_{L^4}\bigr).
\end{equation}
As the assumptions made in the proposition ensure that 
 $\omega_i\in L^\infty([0,T];W^{1,4})$ and that  $u_i\in
L^\infty([0,T];L^2),$  we deduce that
 $\du\in L^\infty([0,T]\times\R^2).$  Therefore, there exists some $T_0>0$ 
such that the right-hand side of the above inequality 
tend to $0$ when  $p$ goes to infinity.  This yields uniqueness
on $[0,T_0].$ {}From a standard connectivity argument, it is now easy
to conclude to uniqueness on the whole interval  $[0,T].$  \end{p}


\subsection{Anisotropic a priori estimates}

If in addition to the $H^1$ hypothesis on $(\theta_0,u_0)$,
 we assume that  $\omega_0\in L^p$ for some  $p$ in 
 $[2,\infty[,$
then the vorticity equation
 \begin{equation}\label{eq:vorticity}
\partial_t\omega+u\cdot\nabla\omega=\partial_1\theta
\end{equation}
implies that
\begin{equation}\label{eq:vor1}
\|\omega(t)\|_{L^p}\leq \|\omega_0\|_{L^p}
+\int_0^t\|\partial_1\theta\|_{L^p}.
\end{equation}
Now, remind that as $\theta_0\in H^1,$ a bound
 for  $\d_1\theta$ in $L^2_{loc}(\R_+;H^1)$ is available, 
whence also in $L^2_{loc}(\R_+;L^p)$
 by Sobolev embedding. 
In fact, we even have a more accurate information 
if  $\omega_0\in\sqrt L.$ 
Indeed,  Lemma  \ref{l:injection} ensures that
 $H^1$ is continuously embedded in  $\sqrt L$
so that, according to \eqref{eq:vor1}, 
\begin{equation}\label{eq:omegaL}\|\omega(t)\|_{\sqrt L}\leq
\|\omega_0\|_{\sqrt L} +C\sqrt{t}\, 
\|\d_1\theta\|_{L^2_t(H^1)}.\end{equation}

However, this bound \emph{does not} imply that   
$\nabla u\in L^1_{loc}(\R_+;L)$
so that one cannot get uniqueness by a direct application of Proposition
\ref{p:uniqueness}. 
In fact, thanks to  \eqref{eq:bs}, it is obvious that 
 $\nabla u\in L^1_{loc}(\R_+;L)$ provided  $\omega\in
L^1_{loc}(\R_+;L^2\cap L^\infty).$
According to \eqref{eq:vor1}, having
  $\d_1\theta$ in $L^1(\R_+;L^\infty)$ will entail that 
the vorticity is bounded.
 \smallbreak
 In order to get this, we shall first show that one may propagate some
additional   horizontal Sobolev regularity for $\theta.$
 By virtue of Lemma \ref{inclusion-sobolev} (see the appendix), this will enable
 us to estimate $\d_1\theta$ in  $L^1(\R_+;L^\infty)$
(and even in  $L^2_{loc}(\R_+;L^\infty)$ actually).

More precisely, we assume from now on that   $(\theta_0,u_0)\in H^1(\R^2)$ and 
$\omega_0\in\sqrt L,$ and that, in addition, 
$|\partial_1|^{1+s}\theta_0\in L^2$ for some  $s\in (0,\frac 12].$
 In order to  propagate the additional regularity,  one
may apply operator $|\partial_1|^{1+s}$  to the equation
$$\partial_t\theta+u\cdot\nabla\theta-\kappa\partial_1^2\theta=0,$$ and take
the $L^2(\R^2)$ inner product with  $|\partial_1|^{1+s}\theta$. After  integrating by
parts, we find that
\begin{equation}\label{eq:aniso}
\frac{1}{2}\frac{d}{dt}\||\partial_1|^{1+s}\theta\|_{L^2}^2
+\kappa||\partial_1|^{2+s}\theta\|_{L^2}^2\leq\big|
(|\partial_1|^{1+s}\big(u\cdot\nabla\theta),|\partial_1|^{1+s}\theta\big)_{L^2}\big|.
\end{equation}
Bounding the right-hand side is the main difficulty.
First of all, let us notice that   $i|\d_1|=\d_1R_1$ where  $R_1$ 
stands for the Riesz operator with respect to the first variable.
As  $|\d_1|^s$ is a symmetric operator, one may write
 $$\big|
(|\partial_1|^{1+s}\big(u\cdot\nabla\theta),|\partial_1|^{1+s}\theta\big)_{L^2}\big|
=\big|
(\partial_1\big(u\cdot\nabla\theta),R_1|\partial_1|^{1+2s}\theta\big)_{L^2}\big|\leq I_1
+I_2
$$
with 
$$I_1:=\big|
\big(\partial_1u\cdot\nabla\theta,R_1|\partial_1|^{1+2s}\theta\big)_{L^2}\big|
\quad\hbox{and}\quad
I_2:=\big|\bigl(u\cdot\nabla\partial_1\theta,R_1|\partial_1|^{1+2s}\theta\bigr)_{L^2}\bigr|.
$$
The term $I_2$ is easy to deal with.
Indeed, for $s\in(0,1/2]$, we have, according to  H\"older and
 Parseval inequalities,
$$
\begin{array}{lll}
I_2&\leq&
\|u\|_{L^\infty}\|\nabla\partial_1\theta\|_{L^2}
\||\partial_1|^{1+2s}\theta\|_{L^2},\\[1ex]
&\leq&
\|u\|_{L^\infty}\|\partial_1\theta\|_{H^1}^2.
\end{array}
$$
Thanks to  \eqref{eq:uLinfty} and by virtue of 
Inequalities \eqref{eq:u1},
\eqref{eq:estimation:h1} and \eqref{eq:omegaL}, we thus have
\begin{equation}\label{eq:I2}
\int_0^tI_2(\tau)\,d\tau\leq C(t,\kappa,\theta_0,u_0)
\end{equation}
where, from now on,   $C(t,\kappa,\theta_0,u_0)$ denotes 
a positive continuous function depending only on $t,$ $\kappa$
and on the norm of  $(\theta_0,u_0)$
in  $H^1\times\bigl(H^1\cap W^{1,4}\bigr).$
\smallbreak
In order to bound the term $I_1$ one may write  $I_1\leq I_1^1+I_1^2$ with  
$$
\begin{array}{lll}
I_1^1&:=&\big|
\big(\partial_1u^1\partial_1\theta,R_1|\partial_1|^{1+2s}\theta\big)_{L^2}\big|,\\[1ex]
I_1^2&:=&\big|\big(\partial_1
u^2\partial_2\theta,R_1|\partial_1|^{1+2s}\theta\big)_{L^2}\big|. \end{array} $$
 For $I_1^1,$ as $\partial_1 u_1=-\partial_2 u_2,$ integrating by parts yields
$$
\begin{array}{lll}
I_1^1&\leq&\bigl|\bigl(u^2\partial_2\partial_1\theta,
R_1|\partial_1|^{1+2s}\theta\bigr)_{L^2}\bigr| +
\bigl|\bigl(u^2\partial_1\theta,
R_1|\partial_1|^{1+2s}\d_2\theta\bigr)_{L^2}\bigr|,\\[1ex] &\leq&
\tilde I_1^1+\hat I_1^1
\end{array}
$$
with
$$\tilde I_1^1:=\bigl|\bigl(u^2\partial_2\partial_1\theta,
R_1|\partial_1|^{1+2s}\theta\bigr)_{L^2}\bigr|\ \hbox{ and }\ 
\hat I_1^1:=\bigl|\bigl(|\partial_1|^{2s}(u^2\partial_1\theta),
R_1|\partial_1|\partial_2\theta\bigr)\bigr|.
$$
In order to bound the  term $\tilde I_1^1,$ one may combine
 H\"older Inequality and \eqref{eq:uLinfty}. As $0<s\leq1/2,$ we get
$$\begin{array}{lll}
\tilde I_1^1&\leq& \|u\|_{L^\infty}\|\partial_1\theta\|_{H^1}
\||\partial_1|^{1+2s}\theta\|_{L^2},\\[1ex]
&\leq&C\bigl(\|u\|_{L^2}+\|\omega\|_{L^4}\bigr)\|\d_1\theta\|_{H^1}^2.
\end{array}
$$
In consequence, by virtue of \eqref{eq:u1},
\eqref{eq:estimation:h1} and \eqref{eq:omegaL}, we have
\begin{equation}\label{eq:I111}
\int_0^t\tilde I_1^1(\tau) d\tau\leq C(t,\kappa, \theta_0,u_0).
\end{equation}
As for $\hat I_1^1,$ we use the fact that 
\begin{equation}\label{eq:I1121}
\hat I_1^1\leq \||\partial_1|^{2s}
(u^2\partial_1\theta)\|_{L^2}\|\partial_1\partial_2\theta\|_{L^2}.
\end{equation}
Because $s\in(0,1/2],$ we have 
$$
\||\partial_1|^{2s}
(u^2\partial_1\theta)\|_{L^2}\leq \|u^2\partial_1\theta\|_{H^1},
$$ whence
$$\begin{array}{lll}
\||\partial_1|^{2s}
(u^2\partial_1\theta)\|_{L^2}&\leq& 
 \|u^2\partial_1\theta\|_{L^2}
+\|\d_1\theta\nabla u^2\|_{L^2}+\|u^2\,\d_1\nabla\theta\|_{L^2},\\[1ex]
&\leq&\|u\|_{L^\infty}\|\d_1\theta\|_{L^2}+\|\nabla u^2\|_{L^4}\|\d_1\theta\|_{L^4}
+\|u\|_{L^\infty}\|\nabla\d_1\theta\|_{L^2}.
\end{array}
$$
Thanks to the Sobolev embedding  $H^1\hookrightarrow L^4$ and to  \eqref{eq:bs},
\eqref{eq:uLinfty}
we get 
$$\||\partial_1|^{2s}
(u^2\partial_1\theta)\|_{L^2}\leq
C\bigl(\|u\|_{L^\infty}+\|\omega\|_{L^4}\bigr)
\|\d_1\theta\|_{H^1}.
$$
Coming back to \eqref{eq:I1121} and using \eqref{eq:uLinfty},
one can now conclude that
\begin{equation}\label{eq:I112}
\int_0^t\widehat I_1^1(\tau)d\tau\leq
C(t,\kappa,\theta_0,u_0).
\end{equation}
The term $I_1^2$ is more intricate to deal with. To start with, 
we integrate by parts to rewrite this term as follows:
$$I_1^2\leq \biggl|\int 
u_2\:\partial_1\partial_2\theta\:R_1|\partial_1|^{1+2s}\theta\,dx\biggr|+\biggl|\int
|\partial_1|^{s}(u_2\partial_2\theta)|\partial_1|^{2+s}\theta\,dx\biggr|,$$ 
from which we get the following bound:
$$I_1^2\leq \|u_2\|_{L^\infty}\|\partial_1\partial_2\theta\|_{L^2}
\||\partial_1|^{1+2s}\theta\|_{L^2}
+\||\partial_1|^{s}(u_2\partial_2\theta)\|_{L^2}\||\partial_1|^{2+s}\theta\|_{L^2}.$$
As $s\in(0,1/2],$ Young inequality enables us to write
\begin{equation}\label{eq:I121}
I_1^2\leq \|u\|_{L^\infty}\|\d_1\theta\|_{H^1}^2
+\frac\kappa2\||\d_1|^{2+s}\theta\|_{L^2}^2
+\frac1{2\kappa}\||\d_1|^s(u_2\d_2\theta)\|_{L^2}^2.
\end{equation}
Let us admit (see the proof in appendix) that there exists a constant $C$ such that
for all $s\in(0,1/2]$ we have  
\begin{equation}\label{eq:anisobound2}
\||\partial_1|^{s}(u_2\d_2\theta)\|_{L^2}
\leq C \|u\|_{H^1}\bigl(\|\d_2\theta\|_{L^2}+\|\d_1\d_2\theta\|_{L^2}\bigr).
\end{equation} 
Using  \eqref{eq:uLinfty} and plugging  
\eqref{eq:u1},
\eqref{eq:estimation:h1} and \eqref{eq:omegaL} in \eqref{eq:I121}, we get
\begin{equation}\label{eq:I12}
\int_0^tI_1^2(\tau)\,d\tau\leq C(t,\kappa,\theta_0,u_0)
+\frac\kappa2\int_0^t\||\d_1|^{2+s}\theta(\tau)\|_{L^2}^2\,d\tau.
\end{equation}
It is now suitable  to integrate  \eqref{eq:aniso} with respect to time
and to plug \eqref{eq:I2}, \eqref{eq:I111}, \eqref{eq:I112} in
\eqref{eq:I12}. We eventually get for all $s\in(0,1/2],$
\begin{equation}\label{eq:estimation:h1+s}
\||\d_1|^{1+s}|\theta(t)\|_{L^2}^2+\kappa\int_0^t\||\d_1|^{2+s}\theta\|\,d\tau
\leq \||\d_1|^{1+s}|\theta_0\|_{L^2}^2+C(t,\kappa,\theta_0,u_0).
\end{equation}
Resorting to Lemma \eqref{inclusion-sobolev} with $s_1=1+s$ and $s_2=1,$
we find that 
$$
\int_0^t\|\d_1\theta\|_{L^\infty}^2\,d\tau
\leq C\int_0^t\Bigl(\|\d_1\theta\|_{L^2}^2+\||\d_1|^{2+s}\theta\|_{L^2}^2
+\|\d_1\d_2\theta\|_{L^2}^2\Bigr)\,d\tau.
$$
Therefore, by virtue of  Inequalities
\eqref{eq:estimation:h1} and \eqref{eq:estimation:h1+s},
we get  a bound for  $\d_1\theta$ in
$L^2([0,t];L^\infty)$ in terms of  $t$ and of the norms of the initial 
data.  As explained before, 
this supplies the desired  bound for the vorticity 
in $L^\infty_{loc}(\R_+;L^\infty).$


\subsection{A global existence result}

This paragraph is devoted to proving the following result (which
obviously implies Theorem \eqref{th:resultat2}):
\begin{theorem}\label{th:main3}
Let $(\theta_0,u_0)\in H^1$ with  $\div u_0.$ 
System $\eqref{eq:b-temperature}$ has a global solution $(\theta,u)$
such that 
$$
(\theta,u)\in\cC_w(\R_+;H^1)\quad\hbox{and}\quad\d_1\theta\in
L^2_{loc}(\R_+;H^1).
$$
If in addition  $\omega_0\in\sqrt L$ then one may construct 
a global solution which also satisfies
$$
\omega\in L^\infty_{loc}(\R_+;\sqrt L).
$$
If in addition  $\omega_0\in L^\infty$ and there exists $s\in(0,1/2]$ such that
$|\d_1|^{1+s}\theta_0\in L^2$ then the above solution is unique, 
strongly continuous in time with values in  $H^1,$ and satisfies 
$$
|\d_1|^{1+s}\theta\in \cC(\R_+;L^2)\quad\hbox{and}\quad
|\d_1|^{2+s}\theta\in L^2_{loc}(\R_+;L^2).
$$
\end{theorem}
\begin{p}
The result may be obtained by means of the 
Friedrichs method. With the  notation of the previous  section, we
solve the following ODE in  $L^2_n$:
\begin{equation*}
\begin{cases}
\partial_t\theta+J_n\div(J_n\theta\,J_nu)-\kappa\d_1^2\cP J_{n}\theta=0,\\
\partial_tu+\cP J_{n}\div(\cP J_{n}u\otimes \cP J_{n}u)
=\cP J_{n}(\theta e_2),\\
(\theta,u)|_{t=0}=J_{n}(\theta_0,u_0).
\end{cases}
\end{equation*}
Cauchy-Lipschitz theorem gives
a unique maximal solution $(\theta_n,u_n)$
in the space $\cC^1([0,T_n^*);L^2_n).$
As  $J_n^2=J_n,$ $\cP^2=\cP$ et
$ J_n\cP=\cP J_n,$ we deduce that
 $\cP u_n=u_n,$ $J_nu_n=u_n$ and
$J_n\theta_n=\theta_n.$ Therefore
$(\theta_n,u_n)$ satisfies
\begin{equation}\label{eq:approche2}
\begin{cases}
\partial_t\theta_n+J_n\div(\theta_{n}u_n)-\kappa\d_1^2\theta_n=0,\\
\partial_tu_{n}+\cP J_{n}\div(u_{n}\otimes u_{n})=\cP J_{n}(\theta_{n} e_2).
\end{cases}
\end{equation}
As usual , because operators $J_n$ 
and  $\cP J_n$ are orthogonal projectors in all the Sobolev 
spaces, all the previous formal a priori estimates pertaining to Sobolev
norms remind true.  More precisely, we still have \eqref{eq:u1}, 
\eqref{eq:omega1} and  \eqref{eq:estimation:h1} so that \begin{itemize}
\item $(\theta_n)_{n\in\N}$ is bounded in
 $L^\infty_{loc}(\R_+;H^1),$
\item $(\d_1\theta_n)_{n\in\N}$ is bounded in  $L^2_{loc}(\R+;H^1),$
\item $(u_n)_{n\in\N}$ is bounded in 
$L^\infty_{loc}(\R_+;H^1).$
\end{itemize}
This is fully enough to pass to the limit (up to extraction)
in System \eqref{eq:approche2} and to get 
the first part of the theorem. 

In order to construct weak  solutions preserving the 
 $\sqrt L$ and the anisotropic regularities, one may 
 smooth out System \eqref{eq:b-temperature}
 by means of an artificial viscosity. 
 More precisely, we first solve the following system for $\epsilon>0$:
 \begin{equation*} \begin{cases}
\partial_t \theta+u\cdot\nabla \theta-\kappa\partial_1^2\theta -\epsilon\Delta\theta=0\\
\partial_t u+u\cdot\nabla u+\nabla\Pi-\epsilon\Delta u=\theta e_2\\
\div u=0
\end{cases}
\end{equation*}
supplemented with smoothed out initial data  $(\theta_0^\e,u_0^\e).$
   
Resorting  again to  the  Friedrichs method that has been 
used in the case $\e=0,$ and noticing  that the cut-off
operator $J_n$ does not modify the Sobolev estimates, 
we get a global  solution $(\theta^\e,u^\e)$  in 
$$
\cC(\R_+;H^1)\cap L^2_{loc}(\R_+;H^2)
$$
satisfying Inequalities  \eqref{eq:omega1} and
\eqref{eq:estimation:h1} uniformly with respect to $\e.$

Actually, using standard methods, one can check that the $H^2$
regularity controls higher Sobolev norms. 
As the initial data are in $H^\infty,$
the solution $(\theta^\e,u^\e)$ thus belongs to 
all the Sobolev spaces, which will enable us to 
make the following computations rigorous.

The $L^p$ estimate over the vorticity may be proved
by multiplying the vorticity  equation 
$$\partial_t\omega^\e+u^\e\cdot
\nabla\omega^\e-\epsilon\Delta\omega^\e=\partial_1\theta^\e$$ by
$|\omega^\e|^{p-2}\omega^\e$, and performing an integration over $\R^2.$
This gives again 
$$\|\omega^\e(t)\|_{L^p}\leq
\|\omega_0\|_{L^p} +\int_0^t\|\partial_1\theta^\e\|_{L^p} \leq
\|\omega_0\|_{L^p}+C\sqrt{pt}\,  \|\d_1\theta^\e\|_{L^2_t(H^1)}.$$ 
It is also clear that all the anisotropic  Sobolev estimates
remain the same, uniformly with respect to $\e.$
Therefore, having $\e$ tend to $0$ yields the end  of the existence
part of Theorem  \ref{th:resultat4}.

Finally, the uniqueness result
is a mere consequence of Proposition \ref{p:uniqueness}.
\end{p}


\section{Appendix}

\subsection{A few inequalities}

Here we prove a few inequalities which have
been used throughout the paper.
\smallbreak\noindent{\it Proof of Inequality \eqref{eq:log}~:}
For proving \eqref{eq:log}, one may split $\nabla u$ into low and high 
frequencies according to the Littlewood-Paley decomposition. 
More precisely, for any  $N\in\N$ one may write
$$
\nabla u=S_N\nabla u+\sum_{q\geq N}\dq\nabla u.
$$
 We thus have
$$
\|\nabla u\|_{L^\infty}\leq\|S_N\nabla u\|_{L^\infty}
+\sum_{q\geq N}\|\dq\nabla u\|_{L^\infty},
$$
whence, using the definition of  $\|\cdot\|_{LL}$
and Bernstein inequalities,
$$
\|\nabla u\|_{L^\infty}\leq(N+1)\|\nabla u\|_{LL}
+C\sum_{q\geq N}2^q\|\dq\nabla u\|_{L^2}.
$$
Given that $\|\dq\nabla u\|_{L^2}=\|\dq\omega\|_{L^2}$
and that $2-s<0,$ we readily get 
$$
\|\nabla u\|_{L^\infty}\leq(N+1)\|\nabla u\|_{LL}
+C2^{N(2-s)}\|\omega\|_{H^{s-1}}.
$$
Now, if  $C\|\omega\|_{H^{s-1}}\leq\|\nabla u\|_{LL}$ then 
taking $N=0$
obviously yields the desired inequality.
Else, one may choose for $N$ the integer part of 
$$
\frac{1}{s-2}
\log_2\biggl(\frac{C\|\omega\|_{H^{s-1}}}{\|\nabla v\|_{LL}}\biggr)
$$
and we still get the desired result.
\hfill\rule{2.1mm}{2.1mm} 
\begin{lemma}\label{l:injection} 
In dimension two, the Sobolev space $H^1$ continuously embeds
in the space $\sqrt L.$
\end{lemma}
\begin{p}
For any  $p\in[2,\infty[$ and $v\in H^1,$ using the Littlewood-Paley 
decomposition and a  Bernstein inequality
enables us to write
$$
\begin{array}{lll}
\|v\|_{L^p}&\leq&\Sum_{q\geq-1}\|\dq v\|_{L^p},\\[1ex]
&\leq& C\Sum_{q\geq-1} 2^{-\frac{2q}p}\,2^q\|\dq v\|_{L^2},\\[1ex]
&\leq&
C\biggl(\Sum_{q\geq-1}2^{-\frac{4q}p}\biggr)^{\frac12}\|v\|_{H^1},\\[2ex]
&\leq& C\sqrt{p-1}\|v\|_{H^1}, \end{array}
$$
whence the desired result.
\end{p}
\smallbreak\noindent{\it Proof of Inequalities \eqref{eq:anisobound}:} 
For stating the first inequality, 
the starting point 
is the following classical one-dimensional 
Gagliardo-Nirenberg inequality:
 \begin{equation}\label{eq:anisobound3}\|f(x_1,\cdot)\|_{L^\infty_{x_2}}\leq
\|f(x_1,\cdot)\|_{L^2_{x_2}}^{1/2}\|\partial_2f(x_1,\cdot)\|_{L^2_{x_2}}^{1/2}.
\end{equation}
Taking the $L^2_{x_1}$ norm of both 
sides and using Cauchy-Schwarz inequality, we get
$$\|f\|_{L^2_{x_1}(L^\infty_{x_2})}\leq C\|f\|_{L^2(\R^2)}^{\frac 12}\|\partial_2 f\|_{L^2(\R^2)}^{\frac 12}.$$
For proving the second inequality, it is only 
a matter of swapping the roles of  variables $x_1$ and $x_2,$ 
and using Minkowski's inequality.
\hfill\rule{2.1mm}{2.1mm} 
\smallbreak\noindent{\it Proof of Inequality \eqref{eq:anisobound1}:} 
{}From \eqref{eq:anisobound3}, we deduce that
$$
\|f\|_{L^\infty}\leq C\|f\|_{L^\infty_{x_1}(L^2_{x_2})}^{1/2}
\|\d_2f\|_{L^\infty_{x_1}(L^2_{x_2})}^{1/2}.
$$
Applying the second inequality of
 \eqref{eq:anisobound} to $f$ and $\d_2f,$ it is now easy to complete
the proof.
\hfill\rule{2.1mm}{2.1mm} 
\smallbreak\noindent{\it Proof of Inequality
\eqref{eq:anisobound2}:} Obviously, it suffices
to state that
 $$\|f g\|_{L^2_{x_2}(H^{1/2}_{x_1})}\leq
C\|f\|_{H^1}\bigl(\|g\|_{L^2}+\|\d_1g\|_{L^2}\bigr).
$$
For proving the above inequality, we first notice
that the standard product laws for one-dimensional
 Sobolev spaces ensure that for all fixed $x_2,$ we have
$$
\|(fg)(\cdot,x_2)\|_{H^{1/2}(\R)}\leq C\|f(\cdot,x_2)\|_{H^{1/2}(\R)}
\|g(\cdot,x_2)\|_{H^{1}(\R)}.
$$
Therefore
$$
\|f g\|_{L^2_{x_2}(H^{1/2}_{x_1})}\leq C\|f\|_{L^\infty_{x_2}(H^{1/2}_{x_1})}
\|g\|_{L^2_{x_2}(H^{1}_{x_1})}.
$$
Because the trace operator on $x_2=cste$ is continuous from
 $H^1(\R^2)$ to $H^{1/2}(\R),$ we get the desired inequality.
\hfill\rule{2.1mm}{2.1mm} 
\smallbreak 
In the last part of the paper, anisotropic
Sobolev norms have been used several times.
Below, we state a sufficient condition under which
anisotropic Sobolev spaces are embedded in the
set of bounded functions.
\begin{lemma} \label{inclusion-sobolev} 
For any couple  $(s_1,s_2)$ of positive real numbers satisfying
 $1/s_1+1/s_2<2$ there exists a constant $C$ such that  
$$\|u\|_{L^\infty}\leq C\bigl(\|u\|_{L^2}+\||\d_1|^{s_1}u\|_{L^2}
+\||\d_2|^{s_2}u\|_{L^2}\bigr).$$ \end{lemma} 
\begin{p} Using Fourier variables, we see that
$$\|\widehat u(\xi)\|_{L^1(\R^2)}^2\leq  \bigg(\int
(1+|\xi_1|^{2s_1}+|\xi_2|^{2s_2}) |\widehat  u(\xi)|^2\,d\xi\bigg)\times
\bigg(\int (1+|\xi_1|^{2s_1}+|\xi_2|^{2s_2})^{-1}\,d\xi \bigg).$$ 
Therefore, it suffices to show that
 $$\int
(1+|\xi_1|^{2s_1}+|\xi_2|^{2s_2})^{-1}\,d\xi<\infty.$$  
If we make the
change of  variable $$\xi_1=(1+|\xi_2|^{2s_2})^{\frac{1}{2s_1}}\zeta_1$$
we get 
$$\int(1+|\xi_1|^{2s_1}+|\xi_2|^{2s_2})^{-1}\,d\xi=\int
(1+|\xi_2|^{2s_2})^{-1+\frac{1}{2s_1}}(1+\zeta_1^{2s_1})^{-1}\,d\zeta_1\,d\xi_2.$$
This integral is finite whenever $s_1>\frac 12$ and
 $s_2(1-\frac{1}{2s_1})>\frac 12$, a condition which is equivalent
to $1/s_1+1/s_2<2.$
\end{p}


\subsection{Losing a priori estimates}

The second part of the appendix is mainly devoted
to the proof of losing a priori estimates
for the following anisotropic Stokes system with convection
\begin{equation}\label{eq:stokespartiel}
\left\{\begin{array}{l}
\d_tw+v\cdot\nabla w-\nu\d_1^2w+\nabla\Pi=f +g\, e_2,\\
\div w=0\end{array}\right.
\end{equation}
in the case where the gradient of the divergence free
vector field is only in $L^1([0,T];LL^{\frac12})$
(where $LL^{\frac12}$ has been  defined in \eqref{eq:LL1/2}).
Remind that those estimates are the key to the proof of uniqueness
in Theorem \ref{th:resultat1}.
Albeit similar results have been proved before in \cite{DANCHI}, 
we also  prove  losing
a priori estimates for ordinary transport equations for the reader convenience.
 \smallbreak
The key to the proof of all those losing a priori estimates
is the following commutator estimate 
(which is also used in the proof of
Inequality \eqref{eq:apriori}). 
\begin{lemma} Let $v$ be a divergence free vector-field
over $\R^2.$ Let $\omega:=\d_1v^2-\d_2v^1.$ There exists
a positive constant $C$ such that for all  $q\geq-1,$ the term
$F_q(v,\rho):=\ov S_{q-1}v\cdot\nabla\Delta_q\rho
-\Delta_q(v\cdot\nabla\rho)$
(with $\ov S_{q-1}$ defined in \eqref{eq:Sq})
satisfies the following estimates~:
\begin{eqnarray}\label{eq:Fqa}
&\|F_q(v,\rho)\|_{L^2}\leq C\|\nabla v\|_{L^\infty}\Sum_{q'\geq
q-4}2^{q-q'}\|\Delta_{q'}\rho\|_{L^2}+\|\rho\|_{L^\infty}
\sum_{|q-q'|\leq4}\|\Delta_{q'}\omega\|_{L^2},\\
\label{eq:Fqb}
&\|F_q(v,\rho)\|_{L^2}\leq C\sqrt{q+2}\,\|\nabla v\|_{LL^{\frac12}}\Sum_{q'}
2^{-|q-q'|}\|\Delta_{q'}\rho\|_{L^2}.
\end{eqnarray}
In the case $\rho=\omega,$ we have in addition
\begin{equation}\label{eq:Fqc}
\|F_q(v,\omega)\|_{L^2}\leq C\|\nabla v\|_{L^\infty}\Sum_{q'\geq
q-4}2^{q-q'}\|\Delta_{q'}\omega\|_{L^2}.
\end{equation}
\end{lemma}
\begin{p}
Decompose $F_q(v,\rho)$ into
$F_q^1(v,\rho)+F_q^2(v,\rho)+F_q^3(v,\rho)+F_q^4(v,\rho)$ 
with
$$
\begin{array}{lll}
F_q^1(v,\rho):=\Sum_{q'\geq-1}[\ov
S_{q'-1}v,\dq]\cdot\nabla\Delta_{q'}\rho,&&
F_q^2(v,\rho):=\Sum_{q'\geq-1}\bigl(\ov S_{q-1}-\ov S_{q'-1}\bigr)v\cdot
\nabla\dq\Delta_{q'}\rho,\\ F_q^3(v,\rho):=-\dq\Bigl(\Sum_{q'\geq1}
S_{q'-1}\d_i\rho\, \Delta_{q'}v^i\Bigr),&&
F_q^4(v,\rho):=-\!\Sum_{q'\geq0} \d_i\dq\biggl(
\Delta_{q'}v^i\Bigl(\Sum_{|\alpha|\leq1}\Delta_{q'\!+\!\alpha}\Bigr)
\rho\biggr).
\end{array}
$$
Let us emphasize that only the term $F_q^1$ 
involves low frequencies of  $v.$ 
Taking advantage of the support properties
of the function $\varphi$ defined
at the beginning of Subsection \ref{ss:regloc},
we notice that the summation
in the definition of  $F_q^1$ may be restricted to those
indices $q'$ such that  $|q'-q|\leq 4.$
Therefore, a standard commutator inequality
(see e.g. \cite{BCD}, Chap. 2) ensures that
\begin{equation}\label{eq:Fq1}
\Vert F_q^1(v,\rho)\Vert_{L^2}\leq C\Sum_{|q'-q|\leq4}
\Vert{\nabla\ov S_{q'-1}v}\Vert_{L^\infty}
\Vert\Delta_{q'}\rho\Vert_{L^2}.
\end{equation}
For  $F_q^2(v,\rho),$ we obtain, according to H\"older and
 Bernstein inequalities, and to the  localization 
properties of  the Littlewood-Paley decomposition, 
\begin{equation}\label{eq:Fq2}
\Vert F_q^2(v,\rho)\Vert_{L^2}\leq
C\Sum_{|q'-q|\leq1}
\Vert{\nabla\check\Delta_q v}\Vert_{L^\infty}\Vert\dq\rho\Vert_{L^2}
\ \hbox{ with }\ \check
\Delta_q:=\sum_{|\alpha|\leq4}\Delta_{q+\alpha}.
\end{equation}
From the definition of operator $S_{q'-1},$  
the  localization properties of operators $\Delta_q$
and  Bernstein inequalities, we get
\begin{equation}\label{eq:Fq3a}
\Vert F_q^3(v,\rho)\Vert_{L^2}\leq\Sum_{q'\leq q+3}2^{q''-q}
\Vert{\Delta_{q''}\rho}\Vert_{L^2}\Vert{\check\Delta_q\nabla v}\Vert_{L^\infty}.
\end{equation}
Notice that one can alternately  get the following inequality~:
\begin{equation}\label{eq:Fq3b}
\Vert F_q^3(v,\rho)\Vert_{L^2}\leq C\|\rho\|_{L^\infty}\sum_{|q'-q|\leq4} 
\|\Delta_{q'}\omega\|_{L^2}.
\end{equation}
Indeed, it is only a matter of using 
that  the sum defining $F_q^3(v,\rho)$
may be restricted to $q'\geq1$ and thus, according
to  Bernstein inequalities and to
$\|\nabla\Delta_{q'}v\|_{L^2}=\|\Delta_{q'}\omega\|_{L^2},$
one may write
$$\begin{array}{lll}
\|S_{q'-1}\d_i\rho\, \Delta_{q'}v^i\|_{L^2}&\leq&
C\|S_{q'-1}\d_i\rho\|_{L^\infty}
2^{-q'}\|\nabla\Delta_{q'}v^i\|_{L^2},\\[1ex]
&\leq& C\|\rho\|_{L^\infty}\|\Delta_{q'}\omega\|_{L^2}.\end{array}
  $$
Finally the term $F_q^4(v,\rho)$ may 
be bounded as follows:
\begin{equation}\label{eq:Fq4}
\Vert F_q^4\Vert_{L^2}\leq\Sum_{q'\geq q-3}\!2^{q-q'}
\Vert\Delta_{q'}\rho\Vert_{L^2}\Vert\nabla\check\Delta_{q'}v\Vert_{L^\infty}.
\end{equation}
Because
$$
\|\nabla\check\Delta_q v\|_{L^\infty}\leq C\|\nabla v\|_{L^\infty}\quad\hbox{and}\quad
\|\nabla\check\Delta_q v\|_{L^\infty}\leq C\sqrt{q\!+\!2}\,\|\nabla v\|_{LL^{\frac12}},
$$
Inequalities \eqref{eq:Fq1} to  \eqref{eq:Fq4} enable
us to get  \eqref{eq:Fqa} and  \eqref{eq:Fqb}. Inequality   \eqref{eq:Fqc}
stems from \eqref{eq:Fqa}.
\end{p}
One can turn to the statement of losing
a priori estimates.
For technical reasons, we adopt 
the framework of Besov spaces $B^\sigma_{2,\infty}.$
As we have
$H^\sigma\hookrightarrow B^\sigma_{2,\infty}$
and  $B^{\sigma}_{2,\infty}\hookrightarrow H^{\sigma'}$
for all $\sigma>\sigma',$ it is of course not difficult to 
rewrite all those estimates in terms of Sobolev norms.

For the transport equation, we shall prove 
the following result (in the spirit of \cite{BCh,DANCHI}).
\begin{proposition}\label{p:transport}
Let $\rho$ satisfy the transport equation 
\begin{equation}\label{eq:transport}
\partial_t\rho+v\cdot\nabla \rho=f
\end{equation}
with initial data $\rho_0\in B^s_{2,\infty}$ and source term
$f\in L^1([0,T];B^s_{2,\infty}).$ Assume in addition 
that ${\rm div}\,v=0$ and  that, for some $V\in L^1([0,T])$
we have 
\begin{equation}\label{eq:loglip}
\sup_{N\geq0}\frac{\|\nabla S_Nv(t)\|_{L^\infty}}{\sqrt{1+N}}
\leq V(t).
\end{equation}
For all $s\in(-1,1),$ there exists a constant $C$  depending 
only on  $s$ such that
for all $\e\in]0,(s+1)/2[$ and $t\in[0,T],$ we have
$$
\|\rho(t)\|_{B^{s-\e}_{2,\infty}}\leq C\exp\biggl(\frac C\e
\biggl(\int_0^tV(\tau)\,d\tau\biggr)^{\!2}\biggr)\biggl(\|\rho_0\|_{B^s_{2,\infty}}
+\|f\|_{L^1_t(B^s_{2,\infty})}\biggr).
$$
\end{proposition}
\begin{p}
Applying  $\Delta_q$ to Equation \eqref{eq:transport}, 
one may write
$$\partial_t\Delta_q\rho+\ov S_{q-1}v\cdot\nabla\Delta_q\rho
=\Delta_q f+F_q(v,\rho)\quad\hbox{with}\quad
F_q(v,\rho):=
\ov S_{q-1}v\cdot\nabla\Delta_q\rho-\Delta_q(v\cdot\nabla\rho).$$
Taking the  $L^2$ inner product of this inequality 
with  $\Delta_q\rho$ and observing that 
$\div \ov S_{q-1}v=0,$ we thus get
\begin{equation}\label{eq:rho1}
\|\Delta_q\rho(t)\|_{L^2}\leq \|\Delta_q\rho_0\|_{L^2}
+\int_0^t\|\Delta_qf\|_{L^2}\,d\tau+ \int_0^t\|F_q(v,\rho)\|_{L^2}\,d\tau.
\end{equation}
From Inequality \eqref{eq:Fqb}, we readily get
for all    $\e\in]0,(s+1)/2[,$
$q\geq-1$ and  $t\in[0,T],$ 
\begin{equation}
\label{eq:rho2}
2^{q(s-\e)}\|F_q(v(t),\rho(t))\|_{L^2}
\leq C \sqrt{q+2}\,V(t)\|\rho(t)\|_{B^{s-\e}_{2,\infty}}
\end{equation} for some constant  $C$ depending
only on $s.$
\smallbreak
Set
 $\eta=\e/\int_0^TV(\tau)\,d\tau$ 
 and $s_t:=s-\eta\int_0^tV(\tau)\,d\tau$ for  $t\in[0,T].$ 
Putting  \eqref{eq:rho1} and  \eqref{eq:rho2} together 
yields
\begin{eqnarray}\label{eq:rho3}
2^{(2\!+\!q)s_t}\|\Delta_q\rho(t)\|_{L^2}\leq 
2^{(2\!+\!q)s}\|\Delta_q\rho_0\|_{L^2}
2^{-\eta(2\!+\!q)\int_0^tV(\tau')\,d\tau'}\hspace{3cm}\nonumber\\
+\Int_0^t2^{(2\!+\!q)s_\tau}\|\Delta_qf(\tau)\|_{L^2}
2^{-\eta(2\!+\!q)\int_\tau^tV(\tau')\,d\tau'}\,d\tau\qquad\qquad\nonumber\\
\hspace{5cm}+C\Int_0^t\sqrt{2\!+\!q}\, V(\tau)
2^{-\eta(2\!+\!q)\int_\tau^tV(\tau')\,d\tau'} 
\|\rho(\tau)\|_{B^{s_\tau}_{2,\infty}}\,d\tau.
\end{eqnarray}
Notice that  if $q$  satisfies 
\begin{equation}\label{eq:limitedloss1}
2+q\geq\frac{4C^2}{\eta^2\log4}
\end{equation}
then the last term may be bounded by
$$
\frac12\sup_{\tau\in[0,t]}
\|\rho(\tau)\|_{B^{s_\tau}_{2,\infty}}
$$
whereas if  $q$ does not satisfy \eqref{eq:limitedloss1}
then it may be bounded by
$$
\frac{2C^2}{\eta\log2}
\int_0^t  V(\tau) 
\|\rho(\tau)\|_{B^{s_\tau}_{2,\infty}}\,d\tau.
$$
So finally, taking the supremum over $q\geq-1$
in \eqref{eq:rho3}
and using the above two inequalities, we get
$$
\sup_{\tau\in[0,t]} \|\rho(\tau)\|_{B^{s_\tau}_{2,\infty}}\leq
2\|\rho_0\|_{B^s_{2,\infty}}
+2\int_0^t  \|f(\tau)\|_{B^{s_\tau}_{2,\infty}}d\tau
+\frac C\eta\int_0^tV(\tau)
\|\rho(\tau)\|_{B^{s_\tau}_{2,\infty}}\,d\tau.
$$
Thanks to  Gronwall lemma, we end up with
$$
\sup_{t\in[0,T]} \|\rho(t)\|_{B^{s_t}_{2,\infty}}\leq
2e^{\frac C\eta\int_0^TV(t)\,dt}
\biggl(\|\rho_0\|_{B^s_{2,\infty}}
+\int_0^T  \|f(t)\|_{B^{s_t}_{2,\infty}}\,dt\biggr),
$$
which entails the desired inequality
given that  $s\geq s_t\geq s-\e$ for all $t\in[0,T].$
\end{p}
A similar result turns 
to be true for System \ref{eq:stokespartiel}.
In addition, owing to 
the anisotropic viscosity, 
we get an extra horizontal smoothing
(which was the key to the proof
of Proposition \ref{p:unicite}.
More precisely, we have:
\begin{proposition}\label{p:stokespartiel}
Let $v$ and $s$ be as in Proposition $\ref{p:transport}.$
Then we have
$$\displaylines{
\|w(t)\|_{B^{s-\e}_{2,\infty}}+\nu^{\frac12}\|\d_1w\|_{L_t^2(B^{s-\e}_{2,\infty})}
\hfill\cr\hfill\leq C(1+\sqrt{\nu t})
\exp\biggl(\frac C\e \biggl(\int_0^tV(\tau)\,d\tau\biggr)^{\!2}\biggr)
\biggl(\|\rho_0\|_{B^s_{2,\infty}} +\|f\|_{L^1_t(B^s_{2,\infty})}
+\nu^{-\frac12}\|g\|_{L^2_t(B^{s-1}_{2,\infty})}\biggr).}
$$
\end{proposition}
\begin{p}
Let us first apply operator  $\dq$ to 
System \eqref{eq:stokespartiel}.
With the notation introduced 
in the proof of Proposition \ref{p:transport},
we have
$$
\d_t\dq w+\ov S_{q-1}v\cdot\nabla\dq w-\nu\d_1^2\dq w+\nabla\dq\Pi=
\dq f +\dq g\: e_2+F_q(v,w)
$$
with  $F_q(v,w)$ satisfying \eqref{eq:rho2}.
\smallbreak
Taking the  $L^2$ inner product and using the
fact that $\div v=\div w=0,$ we see that
\begin{equation}\label{eq:w2}
\frac12\frac d{dt}\|\dq w\|_{L^2}^2
+\nu\|\d_1\dq w\|_{L^2}^2
=\int \dq f\cdot\dq w\,dx
+\int F_q(v,w)\cdot \dq w\,dx
+\int \dq g \dq w^2\,dx.
\end{equation}
Assume that  $q\geq0.$ 
Taking advantage of Parseval equality, one may write
$$\begin{array}{lll}
\Int \dq g \dq w^2\,dx
&=&-\Int(-\Delta)^{-1}\dq g\,\Delta\dq w^2\,dx,\\
&=&-\Int(-\Delta)^{-1}\dq g\,\dq \d_1^2w^2\,dx
-\Int(-\Delta)^{-1}\dq g\,\dq \d_2^2w^2\,dx.
\end{array}
$$
As $\div w=0,$ integrating by parts yields 
$$\begin{array}{lll}
\Int \dq g \dq w^2\,dx
&=&-\Int(-\Delta)^{-1}\dq g\,\dq \d_1^2w^2\,dx
+\Int(-\Delta)^{-1}\dq g\,\dq \d_1\d_2w^1\,dx
,\\&=&
\Int\d_1(-\Delta)^{-1}\dq g\,\dq \d_1w^2\,dx
-\Int\d_2(-\Delta)^{-1}\dq g\,\dq \d_1w^1\,dx.\end{array}
$$
Next, applying Bernstein and  Young inequalities, we 
deduce that
$$
\int \dq g \dq w^2\,dx\leq C2^{-q}\|\dq g\|_{L^2}\|\d_1\dq w\|_{L^2}
\leq\frac\nu2\|\d_1\dq w\|_{L^2}^2+\frac{C}{2\nu}2^{-2q}\|\dq w\|_{L^2}^2.
$$
Then coming back to \eqref{eq:w2} and integrating,
we thus get for all  $q\geq0,$
$$
\displaylines{
\|\dq w\|_{L_t^\infty(L^2)}^2
+\nu\|\d_1\dq w\|_{L_t^2(L^2)}^2
\leq\|\dq w_0\|_{L^2}^2 \hfill\cr\hfill+
2\|\dq f\|_{L^1_t(L^2)}^2
+2\|F_q(v,w)\|_{L^1_t(L^2)}^2
+\frac C\nu 2^{-2q}\|\dq g\|_{L_t^2(L^2)}^2.}
$$
For  $q=-1,$ we merely have
$$
\|\Delta_{-1}w(t)\|_{L^2}
\leq\|\Delta_{-1}w_0\|_{L^2}
+\int_0^t\Bigl(\|\Delta_{-1}f\|_{L^2}
+\|\Delta_{-1}g\|_{L^2}+\|F_{-1}(v,w)\|_{L^2}\Bigr)\,d\tau.
$$
Of course $\|\d_1\Delta_{-1}w\|_{L^2_t(L^2)}
\leq Ct^{\frac12}\|\Delta_{-1}w\|_{L^\infty_t(L^2)}.$
So finally, for all $q\geq-1,$ we have
$$\displaylines{
\|\dq w\|_{L_t^\infty(L^2)}
+\nu^{\frac12}\|\d_1\dq w\|_{L_t^2(L^2)}
\hfill\cr\hfill\leq 2(1+\sqrt{\nu t})
\biggl(\|\dq w_0\|_{L^2}+\|\dq f\|_{L^1_t(L^2)}+\|F_q(v,w)\|_{L^1_t(L^2)}
+\frac C{\nu^{\frac12}} 2^{-q}\|\dq g\|_{L_t^2(L^2)}\biggr).}
$$
With Inequality \eqref{eq:rho2} at our disposal,
it is now easy to conclude the proof
of the proposition.
It is just a matter of 
arguing exactly as in Proposition \ref{p:transport}.
\end{p}


  \end{document}